\documentclass[final,leqno]{siamltex}
\pdfoutput=1

\usepackage{amsmath}
\usepackage{bbm,amssymb}
\usepackage{graphicx}
\usepackage{xy}
\DeclareGraphicsExtensions{.jpg,.pdf,.png,.jpeg}
\newtheorem{thm}{Theorem}[section]
\newtheorem{cor}[thm]{Corollary}
\newtheorem{lem}[thm]{Lemma}

\newtheorem{remark}[thm]{Remark}

\newcommand{\p}{\partial}
\newcommand{\nab}{\nabla}
\newcommand{\Ome}{\Omega}
\newcommand{\Del}{\Delta}

\newcommand{\bn}{\overline{n}}
\newcommand{\mynorm}[1]{\overline{| #1 |}}

\newcommand{\lt}{{L^2}}

\begin{document}

\title{Finite element methods for a bi-wave equation modeling d-wave 
superconductors}

\author{Xiaobing Feng\thanks{Department of Mathematics, The University of
Tennessee, Knoxville, TN 37996, U.S.A. (xfeng@math.utk.edu).}
\and
Michael Neilan\thanks{Department of Mathematics, The University of
Tennessee, Knoxville, TN 37996, U.S.A. (neilan@math.utk.edu).}
}

\markboth{X. Feng and M. Neilan}{Finite element methods for a bi-wave equation}

\maketitle

\begin{abstract}
In this paper we develop two conforming finite 
element methods for a fourth order bi-wave equation arising 
as a simplified Ginzburg-Landau-type model for $d$-wave superconductors
in absence of applied magnetic field. 
Unlike the biharmonic operator $\Delta^2$, the bi-wave 
operator $\Box^2$ is not an elliptic operator, so the energy
space for the bi-wave equation is much larger than the energy
space for the biharmonic equation. This then makes it possible to
construct low order conforming finite elements for the 
bi-wave equation. However, the existence and construction of 
such finite elements strongly depends on the mesh. In the paper,
we first characterize mesh conditions which allow and not 
allow construction of low order conforming finite elements
for approximating the bi-wave equation. We then construct a cubic 
and a quartic conforming finite element. It is proved that both 
elements have the desired approximation properties, and 
give optimal order error estimates in the energy norm,
suboptimal (and optimal in some cases) order error estimates 
in the $H^1$ and $L^2$ norm.
Finally, numerical experiments are presented to guage the
efficiency of the proposed finite element methods and to
validate the theoretical error bounds.
\end{abstract}



\section{Introduction}\label{sec-1}
This paper concerns finite approximations of the following boundary 
value problem:
\begin{alignat}{2}\label{boxproblem1}
\delta \Box^2 u-\Delta u&=f\qquad &&\text{in }\Omega,\\
u=\p_{\bar{n}} u&=0\qquad 
&&\text{on }\partial \Omega, \label{boxproblem2}
\end{alignat}
where $0<\delta\ll 1$ is a given (small) number,
\begin{alignat*}{2}
\Box u &:=\partial_{xx}u-\partial_{yy}u, \qquad\qquad \Box^2u &&:=\Box(\Box u), \\
\bar{n} &:=(n_1,-n_2), \qquad\qquad \p_{\bar{n}} u &&:=\nab u\cdot \bn,
\end{alignat*}
$\Omega\subset \mathbf{R}^2$ is a bounded domain with piecewise smooth
boundary $\p\Ome$, and $n:=(n_1,n_2)$ denotes the unit outward normal 
to $\partial\Omega$. 
As $\Box$ is the well-known ($2$-D) wave operator, we shall call $\Box^2$ 
the bi-wave operator throughout this paper. It is easy to verify that
\[
\Box^2u(x,y)=\frac{\p^2 u}{\p x^2} -2\frac{\p^2 u}{\p x\p y} 
              +\frac{\p^2 u}{\p y^2}.
\]
Hence, equation \eqref{boxproblem1} is a fourth order PDE, which can
be viewed as a singular perturbation of the Poisson equation by the 
bi-wave operator. As a comparison, we recall that the biharmonic 
operator $\Del^2$ is defined as
\[
\Del^2 u(x,y):=\Del(\Del u(x,y))=\frac{\p^2 u}{\p x^2} +2\frac{\p^2 u}{\p x\p y} 
              +\frac{\p^2 u}{\p y^2}.
\]
Although there is only a sign difference in the mixed derivative term,
the difference between $\Del^2$ and $\Box^2$ is 
fundamental because $\Del^2$ is an elliptic operator while 
$\Box^2$ is a hyperbolic operator.

Superconductors are materials that have no resistance to the 
flow of electricity when the surrounding temperature is 
below some critical temperature. At the superconducting state, 
the electrons are believed to ``team up pairwise"
despite the fact that the electrons have negative charges and normally 
repel each other. The Ginzburg-Landau theory \cite{GL50} has been well 
accepted as a good mean field theory for low (critical temperature) $T_c$ 
superconductors \cite{tinkham04}. However, a theory to explain high $T_c$ 
superconductivity still eludes modern physics. 
In spite of the lack of satisfactory microscopic 
theories and models, various generalizations of the Ginzburg-Landau-type
models to account for high $T_c$ properties such as the anisotropy and 
the inhomogeneity have been proposed and developed. 
In low $T_c$ superconductors, electrons are thought to pair 
in a form in which the electrons travel together in spherical orbits, 
but in opposite directions. Such a form of pairing is often called
$s$-wave \cite{tinkham04}. However, in high $T_c$ superconductors, 
experiments have produced strong evidence for $d$-wave pairing
symmetry in which the electrons travel together in orbits resembling 
a four-leaf clover (cf. \cite{Du99,rxt96,xrt96,fk97} and the 
references therein). Recently, the $d$-wave pairing has gained substantial 
support over $s$-wave pairing as the mechanism by which 
high-temperature superconductivity might be explained. 
In generalizing the Ginzburg-Landau models to high $T_c$
superconductors, the key idea is to introduce multiple order 
parameters in the Ginzburg-Landau free energy functional. 
These models, which can also be derived from the phenomenological
Gorkov equations \cite{fk97}, have built a reasonable basis 
upon which detailed studies of the fine vortex structures
in some high $T_c$ materials have become possible. 
We refer the reader to \cite{Du99,rxt96,xrt96,fk97} and the 
references therein for a detailed exposition on modeling and analysis
of $d$-wave superconductors.

We obtain equation \eqref{boxproblem1} from the Ginzburg-Landau-type
$d$-wave model considered in \cite{Du99} (also see 
\cite{rxt96,xrt96}) in absence of applied magnetic field
by neglecting the zeroth order nonlinear terms but retaining the 
leading terms.  In the equation, $u$ (notation $\psi_d$ is instead
used in the cited references) denotes the $d$-wave {\em order parameter}.
We note that the original order parameter $\psi_d$ in the Ginzburg-Landau-type
model \cite{rxt96,Du99} is a complex-valued scalar function whose 
magnitude represents the density of superconducting charge carriers, 
however, to reduce the technicalities and to present the ideas, 
we assume $u$ is a real-valued scalar function in this paper and remark 
that the finite element methods developed in this paper can be easily 
extended to the complex case.  We also note that the parameter $\delta$
appears in the full model as $\delta=-\frac{1}{\beta}$,
where $\beta$ is proportional to the ratio $\frac{\ln(T_{s0}/T)}{\ln(T_{d0}/T)}$
with $T_{s0}$ and $T_{d0}$ being the critical temperatures of
the $s$-wave and $d$-wave components. Clearly, $\beta<0$ (or $\delta>0$)
when $T_{s0}<T<T_{d0}$ and $\beta\searrow -\infty$ (or $\delta\searrow 0$)
as $T\nearrow T_{d0}$. Hence, $\delta$ is expected to be small 
for $d$-wave like superconductors.  

The primary goal of this paper is to develop conforming finite element 
methods for the reduced $d$-wave model \eqref{boxproblem1}. Since the
bi-wave term is the leading term in the full $d$-wave model,
see \cite[Section 4]{Du99}, any good numerical method for 
\eqref{boxproblem1} should be applicable to the full
$d$-wave model. It is easy to see that the energy space for the
bi-wave equation \eqref{boxproblem1} is 
$V:=\{v\in H^1(\Ome);\, \Box v\in L^2(\Ome)\}$ (see Section \ref{prelim}).
Our main task then is to construct finite element subspaces $V^h$
of the energy space $V$ which should be as simple as possible
but also rich enough to have good approximation properties. 
To this end, we note that $H^2(\Ome)\subset V\subset H^1(\Ome)$,
and hence, the desired finite element space $V^h$ should satisfy
$V^h \subset V\subset H^1(\Ome)$. This immediately implies 
that $ V^h \subset C^0(\overline{\Ome})$ (see \cite{Ciarlet78,Brenner}). 
On the other hand, since $V$ is a proper subspace of $H^1(\Ome)$,
the condition $V^h \subset C^0(\overline{\Ome})$ does not
guarantee that $V^h \subset V$.  Hence, $C^0$ (Lagrange) finite element 
spaces are in general not subspaces of $V$.  An intriguing question
is what extra conditions are required to make a $C^0$ finite element
space to be a subspace of $V$. To answer this question, on noting
that $H^2(\Ome)\subset V$, one may choose $V^h$ such that $V^h\subset H^2(\Ome)$,
that is, $V^h$ is a $C^1$ finite element space 
such as Argyris finite element space (cf. \cite[Chapter 6]{Ciarlet78}).
Trivially, $V^h \subset H^2(\Ome) \subset V$. It turns out (see
Section \ref{sec-3}) such a choice would work since it can be shown 
that the finite element solution so defined converges with optimal rate 
in the energy norm of $V$. However, since $C^1$ finite elements 
require either the use of fifth or higher order polynomials with up to
second order derivatives as degrees of freedom \cite{zenisek73,zenisek74}, 
or the use of exotic elements \cite[Chapter 6]{Ciarlet78}, it is
expensive and less efficient to solve the bi-wave equation \eqref{boxproblem1} 
using $C^1$ finite elements. This then motivates us to construct 
low order non-$C^1$ finite elements which give genuine subspaces of $V$
and to develop other types of finite element methods such as 
nonconforming and discontinuous Galerkin methods \cite{feng_neilan}. 
  
The remainder of the paper is organized as follows. Section \ref{prelim}
contains some preliminaries and the functional setting for the 
bi-wave problem. Well-posedness of the problem and 
regularity estimates of the weak solution are established. 
Because $\Box^2$ is a hyperbolic operator, the usual regularity shift 
for fourth order elliptic problems does not hold 
for the bi-wave problem, instead, a weaker shifting 
``rule" only holds. Section \ref{sec-2} devotes to construction
and analysis of piecewise polynomial subspaces of $V$. 
First, we give a characterization of such subspaces. 
It is proved that a subspace of $V$ is ``necessarily" a $C^1$ finite
element space on a general mesh. However, non-$C^1$ finite elements are
possible on restricted meshes. Second, we construct two such 
finite elements. The first one is a cubic element and the second
is a quartic element. Third, we establish the approximation  
properties for both proposed finite elements. Because both elements
are not affine families, a technique of using affine relatives
(cf.  \cite{Brenner,Ciarlet78}) is used to carry out the analysis.
Finally, optimal order error estimates in the energy norm of $V$
are proved for the finite element approximations of  
problem \eqref{boxproblem1}--\eqref{boxproblem2} using the
proposed finite elements. Suboptimal (and optimal in some cases)
order error estimates in the $L^2$-norm are also derived using 
a duality argument. 
In Section \ref{sec-3} we present some numerical
experiment results to gauge the efficiency of the proposed finite element
methods and also to validate our theoretical error bounds.

\section{Preliminaries and functional setting} \label{prelim}
Standard space notation is adopted in this paper.  We refer the reader
to \cite{Brenner,Ciarlet78} for their exact definitions. In addition,
$(\cdot,\cdot)$ and $\langle\cdot, \cdot\rangle_{\partial\Omega}$ are used to denote 
the $L^2$-inner 
products on $\Ome$ and on $\p\Ome$, respectively. $C$ denotes 
a generic $h$ and $\delta$-independent positive constant.  
We also introduce the following special space notation:
\begin{alignat*}{2}
&V_0 :=\{v\in V\cap H^1_0(\Omega);\ \p_{\bar{n}} v\big|_{\partial\Omega}=0\},
&&\qquad
(v,w)_V :=\delta (\Box u, \Box w) +(v,w), \\
&\|v\|_V := \sqrt{(v,v)_V}. &&
\end{alignat*}
It is easy to verify that $(\cdot,\cdot)_V$ is an inner product on $V$, 
hence, $\|\cdot\|_V$ is the induced norm, and $V$ endowed with this inner 
product is a Hilbert space.  We remark that all above claims do not hold
in general if the harmonic term $\Del u$ is dropped in \eqref{boxproblem1} 
because the kernels of the bi-wave operator $\Box^2$ and the wave operator
$\Box$ may contain non-zero functions satisfying the homogeneous Dirichlet 
boundary condition \cite{bd39}.

The variational formulation of \eqref{boxproblem1}--\eqref{boxproblem2} 
can be derived easily by testing \eqref{boxproblem1} against a test
function $v\in V_0$ and using integration by parts formulas. Specifically, 
it is defined as seeking $u\in V_0$ such that
\begin{align}\label{firstvariational}
A^\delta (u,v)=\langle f,v\rangle ,
\end{align}
where 
\[
A^\delta (u,v):= (u,v)_V,
\] 
and $\langle \cdot,\cdot\rangle$ denotes the pairing 
between $V$ and its dual, $V^*$.

We now show that problem \eqref{firstvariational} is well-posed.

\begin{thm}\label{thm2.1}
For any $f\in V^*$, there exists a unique 
solution to \eqref{firstvariational}.  Furthermore, there holds estimate
\begin{align}\label{firstbound}
\|u\|_V\le \|f\|_{V^*}.
\end{align}
\end{thm}

\begin{proof}
We note for $v,w\in V_0$,
\begin{align}\label{coercive}
A^\delta(v,v) &\geq \|v\|_V^2,\\
|A^\delta (v,w)\big| &\leq \|v\|_V\|w\|_V. \label{bounded}
\end{align}  
Then, existence and uniqueness follows directly from an application
of the Lax-Milgram Theorem (cf. \cite{Brenner,Ciarlet78})
and using the fact that $V_0$ is a Hilbert space with the inner 
product $(\cdot,\cdot)_V$.  The estimate \eqref{firstbound} 
follows from \eqref{coercive} and \eqref{firstvariational} after setting
$v=u$ and $w=u$.
\end{proof}

We note $H^2(\Ome)$ is a proper subspace of $V$, so in general 
$u\not\in H^2(\Ome)$ if $f\in V^*$.  However, for smoother function $f$ 
we have the following regularity results.

\begin{thm}\label{regularity_thm}
Assume that the boundary $\p\Ome$ of the domain $\Ome$ is
sufficiently smooth. Let $s_1,s_2$ be two nonnegative integers. Then
there exist constants $C_{s_1,s_2}, \hat{C}_{s_1,s_2}>0$ such that 
the weak solution $u$ of \eqref{firstvariational} satisfies 
\begin{alignat}{2} \label{regularity}
\|\p^{s_1}_x \p^{s_2}_y u\|_{V}
\le &C_{s_1,s_2}\|\p^{s_1}_x \p^{s_2}_y f\|_{V^*}
&&\qquad\mbox{if } \p^{s_1}_x \p^{s_2}_y f\in V^*, \\
\sqrt{\delta} \|\Box^2 \p^{s_1}_x \p^{s_2}_y u\|_{L^2}
+\sqrt{\delta} \|\nab \Box \p^{s_1}_x \p^{s_2}_y u\|_{L^2}
&\,  &&  \nonumber \\
+\|\Del \p^{s_1}_x \p^{s_2}_y u\|_{L^2} 
\le& \hat{C}_{s_1,s_2}\|\p^{s_1}_x \p^{s_2}_y f\|_{L^2}
&&\qquad\mbox{if } \p^{s_1}_x \p^{s_2}_y f\in L^2(\Ome). \label{regularity2}
\end{alignat}
\end{thm}

\begin{proof}
First, we consider the case that $u$ and $f$ have compact support.
Let $w:=\p^{s_1}_x \p^{s_2}_y u$ and $g:=\p^{s_1}_x \p^{s_2}_y f$. 
Because equation \eqref{boxproblem1} is a linear equation, differentiating 
the equation immediately verifies that $w$ and $g$ satisfy
\begin{equation}\label{e2.6}
\delta \Box^2 w -\Del w = g,
\end{equation}
that is, $w$ is a solution of the bi-wave equation with the source
term $g$. Since $u$ is assumed to have a compact support, then 
$w$ also satisfies the homogeneous boundary conditions in
\eqref{boxproblem2}. Thus, it follows from Theorem \ref{thm2.1} that
\[
\|w\|_V\le \|g\|_{V^*},
\]
which gives \eqref{regularity} with $C_{s_1,s_2}=1$.

To show \eqref{regularity2}, it suffices to prove that 
\begin{equation}\label{e2.7}
\sqrt{\delta} \|\nab \Box w\|_{L^2}
+\|\Del w\|_{L^2} \le \hat{C}_{s_1,s_2}\|g\|_{L^2},
\end{equation}
which is equivalent to prove that \eqref{regularity2} holds
for $s_1=s_2=0$. To this end, testing \eqref{e2.6} with $-\Del w$ yields
\[
-\delta(\Box^2 w, \Del w) + \|\Del w\|_{L^2}^2 = -(g, \Del w).
\]
Using the following integral identity
\[
(\Box^2 \varphi, \psi)_\Ome
= \langle \p_{\bar{n}} \Box \varphi, \psi\rangle_{\p\Ome}
  -\langle \Box \varphi, \p_{\bar{n}}\psi\rangle_{\p\Ome}
  + (\Box \varphi, \Box \psi)_\Ome
\]
followed by using Green's identity (for $\Del$)
in the first term on the left hand side we get
\[
-\delta(\Box^2 w, \Del w) = \delta \|\nab \Box w\|_{L^2}^2.
\]
Here, we have dropped the boundary integral terms because $w$
has a compact support.

Combining the above two identities for $w$ and using Schwarz inequality 
yield
\[
\delta \|\nab \Box w\|_{L^2}^2 + \|\Del w\|_{L^2}^2 
\leq \frac12 \|g\|_{L^2}^2 +\frac12 \|\Del w\|_{L^2}^2.
\]
Hence, the above inequality and \eqref{e2.6} imply
that \eqref{e2.7} holds with $\hat{C}_{s_1,s_2}=2\sqrt{2} +1$.

Second, in the case $u$ and $f$ do not have compact support, 
it is clear that $w$ and $g$ still satisfy \eqref{e2.6}.
However, $w$ and its derivatives may not satisfy the homogeneous 
boundary conditions in \eqref{boxproblem2}. To get around this 
difficulty, the well-known tricks are to use the cutoff function technique 
(see \cite{Evans98,GT01}) for interior estimates and to use the 
flattening boundary technique for boundary estimates. The cutoff function 
technique involves testing \eqref{e2.6} by $w\xi$ and $-\Del w\xi$, 
instead of $w$ and $-\Del w$, for a smooth cutoff function $\xi$. 
Integrating by parts on the left hand side and using Schwarz 
inequality and the properties of the cutoff function then yield the
desired interior estimate similar to \eqref{regularity} 
and \eqref{regularity2}. The flattening boundary technique involves
locally mapping the curved boundary into a flat boundary by a 
smooth map (this requires the smoothness of the boundary $\p\Ome$). 
After the desired boundary estimates are obtained in the new 
coordinates, they are then transferred to the solution $w$ in 
the original coordinates. We omit the technical derivations
and refer the interested reader to \cite{Evans98,GT01}
for a detailed exposition of these techniques applying to 
other linear PDEs.
\end{proof}

\section{Construction and analysis of finite element methods}\label{sec-2}

\subsection{Characterization of finite element subspaces of $V$} 
\label{sec-2.1}
Let $\mathcal{T}_h$ be a quasi-uniform triangulation of $\Omega$ 
with mesh size $h\in (0,1)$, and for a fixed $T\in \mathcal{T}_h$, 
let $(\lambda_1^T,\lambda_2^T,\lambda_3^T)$ denote the
barycentric coordinates, and $a_i\ (1\le i\le 3)$ 
denote the vertices of $T$.  We also let $e_i\ (1\le i\le 3)$
denote the edge of $T$ of which $a_i$ is not a vertex, 
and $b_i$ denote the midpoint of edge $e_i$. Define the 
interior and boundary edge sets of $\mathcal{T}_h$ 
\begin{align*}
\mathcal{E}_h^I:&=\{e;\ e\cap \partial\Omega=\emptyset\},\qquad
\mathcal{E}_h^B:=\{e;\ e\cap \partial\Omega\neq \emptyset\}.
\end{align*}
We also set \[\mathcal{E}_h:=\mathcal{E}_h^I\cup \mathcal{E}_h^B,\]
and for $T\in\mathcal{T}_h$,
\begin{align*}
\omega(T):=\mbox{closure}\left(\bigcup_{\p T^\prime\cap \p T\neq \emptyset} 
T^\prime \right).
\end{align*}
For any $e\in \mathcal{E}_h^I$ such that $e=T_1\cap T_2$,
and $v\in H^1(T_1)\cap H^1(T_2)$, define the jumps of
$v$ across $e$ as (assuming the global label of $T_1$ is bigger
than that of $T_2$)
\[
[v]\big|_e :=v^{T_1}\big|_e-v^{T_2}\big|_e, 
\]
where $v^{T_i}=v\big|_{T_i}$, and $[v]\big|_e :=v^{T_1}\big|_e$ if 
$e\in \mathcal{E}_h^B$.

Similarly, for $v\in H^2(T_1)\cap H^2(T_2),\ \alpha\in \mathbf{R}^2$, 
we define the jumps of $\p_\alpha v:=\nabla v\cdot \alpha$ as follows:
\begin{alignat*}{2}
[\p_\alpha v]\big|_e&:=\p_\alpha v^{T_1}\big|_e-\p_\alpha v^{T_2}\big|_e \qquad 
&&e=\partial T_1\cap \partial T_2\in \mathcal{E}^I_h,\\
[\p_\alpha v]\big|_e&:=\p_\alpha v^{T_1}\big|_e \qquad 
&& e=\partial T_1\cap \partial\Omega\in \mathcal{E}_h^B.
\end{alignat*}
We also define the shorthand notation
\[
\overline{\nabla}v:=(v_x,-v_y),\qquad 
\mynorm{\nabla v}:=\nabla v\cdot\overline{\nabla} v.
\] 

In the rest of the paper, we shall often encounter the following 
characterization of the meshes.

\begin{definition}\label{edgetype}
For $e\in\mathcal{E}_h$, let $n$ and $\tau$ denote the outward unit normal
and unit tangent vector of $e$, respectively.  We say that $e$ is a 
{\em type I} edge if 
\begin{align}\label{rect}
\overline{n}=\tau\quad\text{or}\quad \overline{n}=-\tau.
\end{align}
Otherwise, $e$ is called a {\em type II} edge if condition \eqref{rect} 
does not hold.
\end{definition}

\begin{figure}[htb]
\centerline{
\includegraphics[scale=0.3]{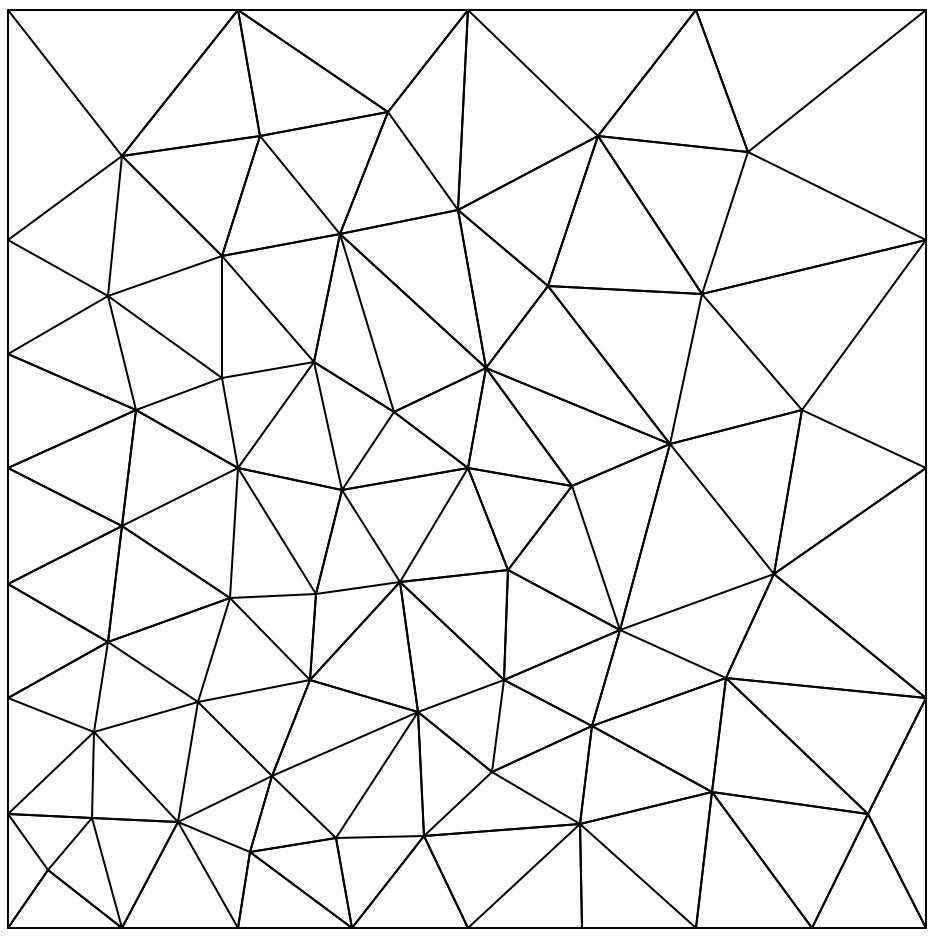}
\qquad\quad
\includegraphics[scale=0.3]{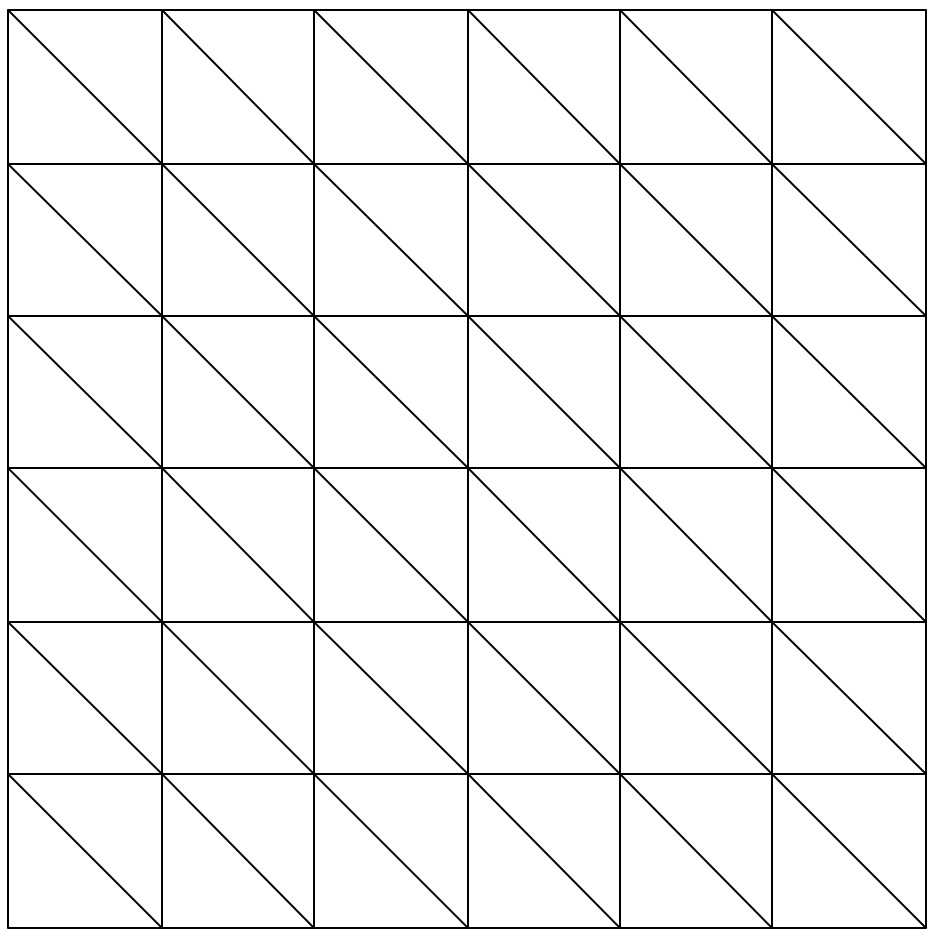}
}
\caption{\label{squaremesh}
{\scriptsize 
Example of meshes of the domain $\Omega=(0,1)^2$ such that every triangle has
no type I edges (left), and one type I edge (right).
}}
\end{figure}

\begin{figure}[htb]
\centerline{
\includegraphics[scale=0.3]{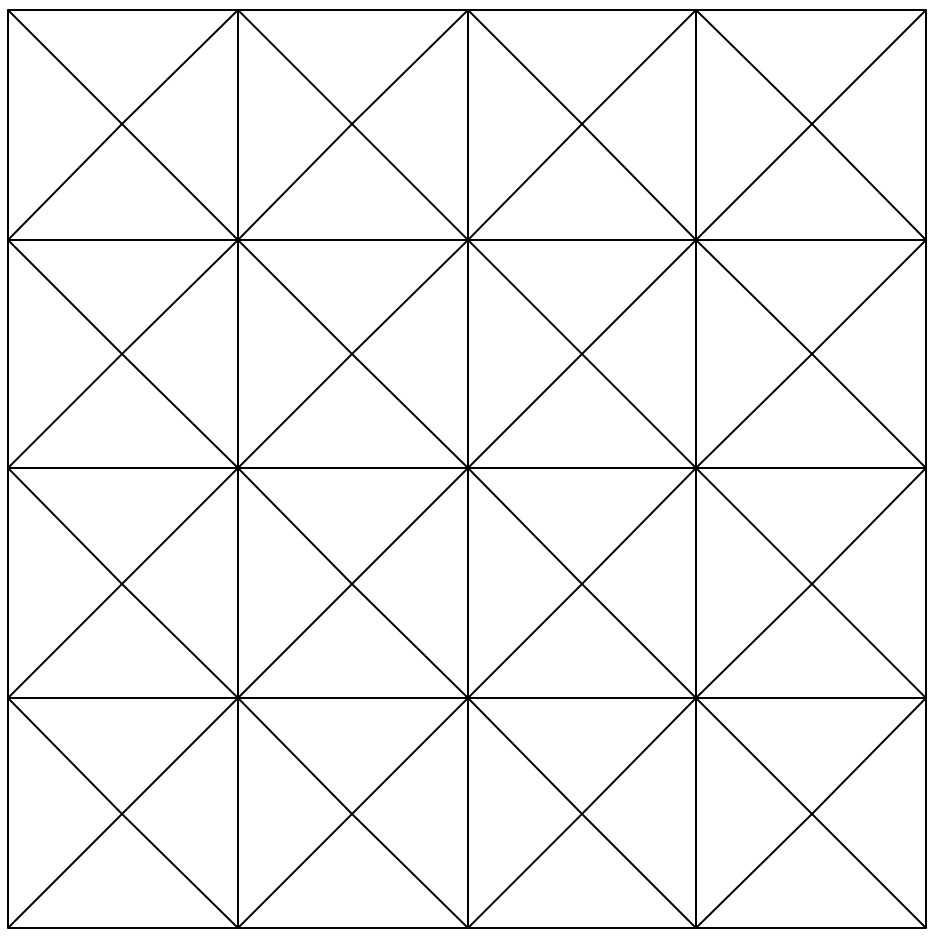}
\qquad\quad
\includegraphics[scale=0.3]{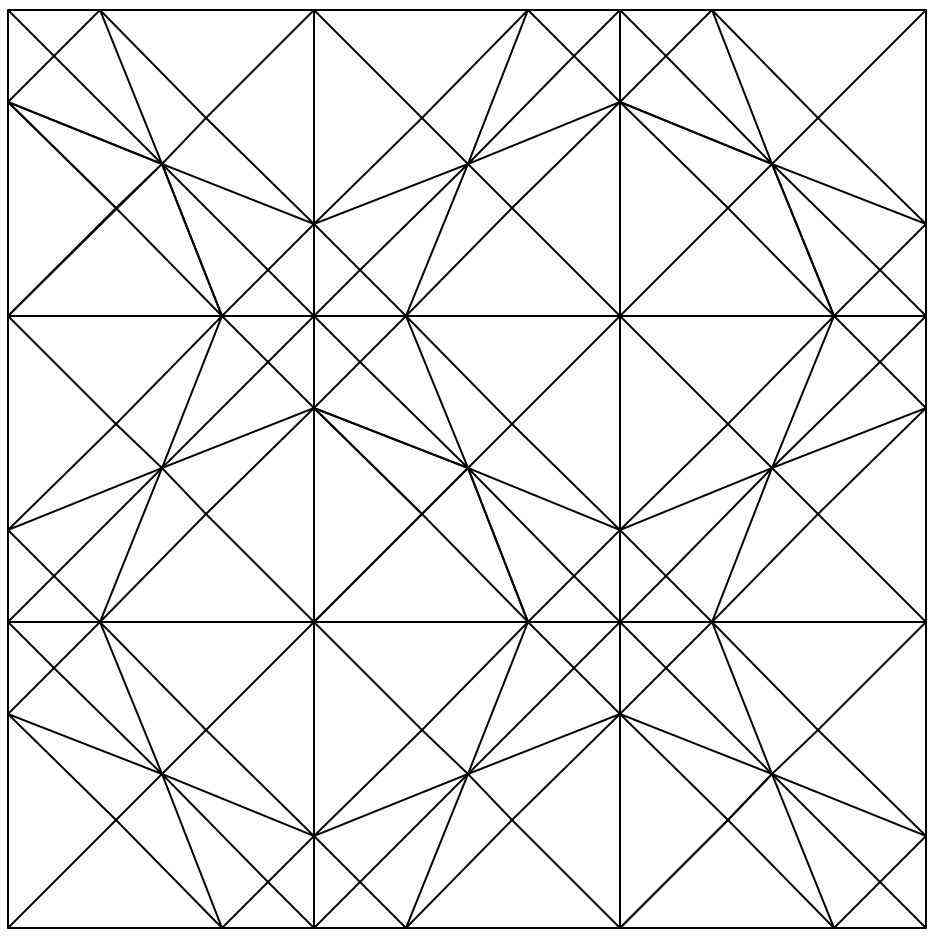}
}
\caption{\label{squaremesh2}{\scriptsize 
Example of a uniform mesh (left) and a nonuniform mesh (right) of 
the domain $\Omega=(0,1)^2$ such that every triangle has two type I edges.}}
\end{figure}

\begin{remark}\label{edgeremark}
(a)  If $e$ is a type I edge, then 
$\overline{n}=(n_1,-n_2)=\pm \tau=\pm (\tau_1,\tau_2)=\pm(n_2,-n_1)$. Therefore,
\[
\tau=\frac{\sqrt{2}}{2}(\pm 1,\pm 1).
\]
That is, the edge $e$ makes an angle of $\frac{\pi}4$ in the plane 
with respect to the $x$-axis. Examples of meshes such that every 
triangle in the partition has exactly zero and one type I edges 
are shown in Figure \ref{squaremesh}, and examples of meshes
such that every triangle has exactly two type I edges are shown 
in Figure \ref{squaremesh2}.

\noindent
(b) For $T\in\mathcal{T}_h$, $e_i\subset \partial T$, 
let $n^{(i)}$ and $\tau^{(i)}$ denote the outward (from $T$) unit normal 
and unit tangent vector of $e_i$, respectively.  Then using the formula
\[
n^{(i)}=-\frac{\nabla \lambda^T_i}{\|\nabla \lambda^T_i\|},
\]
we conclude that $e_i$ is a type I edge if and only if 
\[
\nabla \lambda^T_i\cdot \overline{\nabla}\lambda^T_i=0.
\]
\end{remark}

To construct finite element subspaces of $V$, we first provide the following 
two lemmas, which characterize such spaces.

\begin{lem}\label{H2lem}
Let $X^h$ be a subspace of $V$ consisting of piecewise polynomials, and 
suppose there exists a type II edge $e\in\mathcal{E}_h^I$ with 
$e=\partial T_1\cap \partial T_2$. Then for $v\in X^h$, there holds 
the inclusion $v\in H^2(T_1\cup T_2)$.
\end{lem}

\begin{proof}
Since $X^h$ is finite-dimensional with $X^h\subset H^1(\Omega)$,
we have the inclusion $X^h \subset C^0(\overline{\Ome})$. We also note that
it suffices to show $v \in C^1(\overline{T}_1\cup \overline{T}_2)$ 
for any $v\in X_h$, which in turn is equivalent to show
\[
\big[\p_\alpha v\big]\big|_{e}=0\qquad \forall\alpha\in\mathbf{R}^2.
\]

Let $n$ and $\tau$ denote the normal and tangential direction 
of $e$, respectively. Rewriting $V$ as
\[
V=\{v\in H^1(\Omega);\ \overline{\nabla}v\in H({\rm div};\Omega)\},
\]
there holds for $v\in X^h$
\[
\big[\p_{\bar{n}} v\big]\big|_{e}=0.
\]

Next, using the assumption $\overline{n}\neq \pm \tau$, we can write for any 
constant vector $\alpha\in \mathbf{R}^2$
\begin{align*}
&\big[\p_\alpha v\big]\big|_{e}=\frac{1}{1-(\tau\cdot \overline{n})^2}
\Big\{ \alpha\cdot \big(\tau-\overline{n}(\tau\cdot\overline{n})\big)
\big[\p_\tau v]\Big|_{e} + \alpha\cdot \big(\overline{n}
-\tau(\tau\cdot\overline{n})\big)\big[\p_{\bar{n}} v\big]\Big|_{e}\Big\}.
\end{align*}
But $\big[\p_\tau v\big]\big|_{e}=0$ since $v\in C^0(\overline{\Ome})$ 
and $v\big|_e$ is a polynomial of one variable.  Hence,
$\big[\p_{\bar{n}} v\big]\big|_{e}=0$ implies that
$\big[\p_\alpha v\big]\big|_{e}=0$. The proof is complete.
\end{proof}

\begin{cor}\label{H2cor}
Suppose $X^h$ is a subspace of $V$ consisting of piecewise polynomials, and 
suppose there exists no type I edges in the set $\mathcal{E}_h^I$. Then 
$ X^h\subset H^2(\Omega)$.
\end{cor}

\begin{lem}\label{k5lem}
Suppose $\Sigma_T$ is a linearly independent set of parameters uniquely 
determining a $k$th-degree polynomial $v$ on an interior triangle 
$T\in\mathcal{T}_h$ that includes only function and derivative 
degrees of freedom. Suppose further that $v$ is continuous in 
$\omega(T)$, $\Box v\in L^2(\omega(T))$, and $T$ has at least two 
type II edges that are in the set $\mathcal{E}_h^I$. Then $k\ge 5$.
\end{lem} 

\begin{proof}
If $T$ has three type II edges, then by Lemma \ref{H2lem},
$v\in H^2(\omega(T))$, and it follows that $k\ge 5$ 
(cf. \cite[p.108]{Ciarlet78}, also see \cite{zenisek73,zenisek74}).

Suppose $T$ has exactly two type II edges, without loss of 
generality, assume $e_1$ is type I.  By the proof of 
Lemma \ref{H2lem}, $v$ is $C^1$ across edges $e_2$ and $e_3$.
Let $\mu_i$ denote the order of prescribed derivatives 
at vertex $a_i$ in the set $\Sigma_T$, let $m_i$ denote the 
number of function value (or equivalent) degrees of freedom 
in the set $\Sigma_T$ on edge $e_i$, and let $s_i$ denote the number
of (non-tangential) directional derivative value (or equivalent) 
degrees of freedom in the set $\Sigma_T$ on edge $e_i$.
Since $v$ is continuous in $\omega(T)$, we have
\begin{align}\label{cont}
\mu_2+\mu_3+m_1&\ge k-1,\\
\nonumber\mu_1+\mu_3+m_2&\ge k-1,\\
\nonumber\mu_1+\mu_2+m_3&\ge k-1,
\end{align}
and since $\nabla v$ is continuous across $e_2$ and $e_3$,
\begin{align}
\label{cont2}\mu_1+\mu_2+s_3&\ge k,\\
\nonumber\mu_1+\mu_3+s_2&\ge k.
\end{align}
Adding up the above five inequalities yields
\begin{align*}
4\mu_1+3\mu_2+3\mu_3+m_1+m_2+m_3+s_1+s_2\ge 5k-3.
\end{align*}

Because the set $\Sigma_T$ is linearly independent, and the dimension
of $\Sigma_T$ equals $\frac{(k+1)(k+2)}{2}$, there holds
\begin{align}\label{numdof}
\frac{(k+1)(k+2)}{2}
&\ge \sum_{i=1}^3\Big\{ \frac12(\mu_i+1)(\mu_i+2)+m_i\Big\}+s_2+s_3\\
&\ge \sum_{i=1}^3 \frac12 (\mu_i+1)(\mu_i+2)+5k-3-4\mu_1-3\mu_2-3\mu_3.
\nonumber
\end{align}
Thus,
\begin{align}\label{addline}
(k^2-7k+8)&\ge (\mu_1^2-5\mu_1+2)+(\mu_2-2)(\mu_2-1)+(\mu_3-2)(\mu_3-1).
\end{align}

It is clear that $k$ must be greater than two, therefore, it suffices to show 
that $k$ cannot equal three or four.

\medskip
{\em Case $k=3$:} If $k=3$, by \eqref{addline} we get
\begin{align*}
(\mu_1-3)(\mu_1-2)+(\mu_2-2)(\mu_2-1)+(\mu_3-2)(\mu_3-1)\le 0,
\end{align*}
and since $\mu_i$ are integer-valued, we have
\begin{align*}
1\le \mu_3\le 2,\quad 1\le \mu_2\le 2,\quad 2\le \mu_1\le 3.
\end{align*}
But by \eqref{numdof}, we immediately obtain
\[
10=\frac{(k+1)(k+2)}{2}\ge \sum_{i=1}^3 \frac12 (\mu_i+1)(\mu_i+2)\ge 12,
\]
which is a contradiction.
  
\medskip
{\em Case $k=4$:} As in the previous case, if $k=4$ we have
\begin{align*}
(\mu_1-3)(\mu_1-2)+(\mu_2-2)(\mu_2-1)+(\mu_3-2)(\mu_3-1)\le 0.
\end{align*}
Since
\begin{align*}
1\le \mu_3\le 2,\quad 1\le \mu_2\le 2,\quad 2\le \mu_1\le 3,
\end{align*}
and
\[
15=\frac{(k+1)(k+2)}{2}\ge \sum_{i=1}^3 \frac12 (\mu_i+1)(\mu_i+2),
\]
it is not hard to check that there can only be the following three subcases:  
\begin{align}\label{subcases}
(\mu_1,\mu_2,\mu_3)=(2,1,2),\quad 
(\mu_1,\mu_2,\mu_3)=(2,2,1),\quad
(\mu_1,\mu_2,\mu_3)=(2,1,1).
\end{align} 

If the first subcase holds, then all degrees of freedom lie on the 
vertices, therefore, $m_i,s_i=0,\ 1\le i\le 3$. However, 
it follows from \eqref{cont2} that
\begin{align*}
3=\mu_1+\mu_2\ge 4,
\end{align*}
which is a contradiction. 

A similar argument can be used to exclude the second subcase in \eqref{subcases}. 
Now, suppose $\mu_1=2,\ \mu_2=1,$ and $\mu_3=1$.  By \eqref{cont} 
and \eqref{cont2}, we have
\begin{align*}
m_3\ge 1,\qquad s_3\ge 2,\qquad s_2\ge 1.
\end{align*}
But this implies that
\begin{align*}
15\ge \sum_{i=1}^3 \big\{\frac12 (\mu_i+1)(\mu_i+2)+m_i\big\} +s_2+s_3
\ge 16,
\end{align*}
a contradiction. Thus, the third subcase can not happen, either.
Therefore, we must have $k\ge 5$. The proof is complete.
\end{proof}

By Lemmas \ref{H2lem} and \ref{k5lem}, and Corollary \ref{H2cor}, 
we conclude that unless certain types of meshes are used, we must resort 
to either $C^1$ finite elements such as Argyris, Hsieh-Clough-Tocher,
Bogner-Fox-Schmit elements (cf. \cite{Brenner,Ciarlet78}), or special 
exotic elements (e.g. macro elements), or nonconforming elements
(cf. \cite{feng_neilan}) to solve problem 
\eqref{boxproblem1}--\eqref{boxproblem2},  However for special meshes, 
we now show in the following subsections that it is feasible
to construct low order finite element subspaces of $V$.

\subsection{A cubic conforming finite element}
To construct a cubic conforming finite element, we assume 
that $\mathcal{T}_h$ is a triangulation of $\Ome$ and every triangle 
of $\mathcal{T}_h$ has {\em two} type I edges. Examples of 
such meshes are shown on a square domain in Figure \ref{squaremesh2}.
Our cubic finite element $S^h_3:=(T,P_T,\Sigma_T)$ is defined as follows:
\begin{enumerate}
\item[\rm (i)] $T$ is a triangle with two type I edges,
\item[\rm (ii)] $P_T=\mathbb{P}_3(T)$, the space of cubic polynomials on $T$,
\item[\rm (iii)] 
$\displaystyle{
\Sigma_T=\left\{\begin{array}{ll} v(a_i) &1\le i\le 3,\\
v(a_{i3}) & 1\le i \le 2,\\
\nabla v(a_i)\cdot(a_j-a_i)&1\le i\le 2,\ 1\le j\le 3,\ j\neq i,\\
\p_{\bar{n}} v(b_{3}),
\end{array}\right.
}$

where $e_3$ is a type II edge.
\end{enumerate}

\begin{figure}[htb]
\begin{center}
\indent \indent\includegraphics[angle=0,width=6.15cm,height=4.75cm]{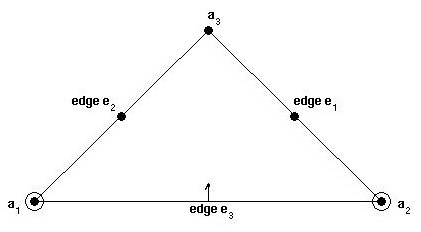}
\end{center}
\caption{\scriptsize Element $S^h_3$. Solid dots indicate function 
evaluation, circles indicate first derivative evaluation, and arrows 
indicate evaluation of derivatives in the direction $\overline{n}$.}
\label{cubic_element}
\end{figure}

\begin{lem}\label{UnisolventLem1}
The set $\Sigma_T$ is unisolvent.  That is, any polynomial of degree 
three is uniquely determined by the degrees of freedom in $\Sigma_T$.
\end{lem}

\begin{proof}
Suppose $v\in\mathbb{P}_3(T)$ equals zero at all the degrees of 
freedom in $\Sigma_T$. To complete the proof, it suffices to show $v\equiv 0$
since dim($\mathbb{P}_3(T))=\text{dim}(\Sigma_T)=10$.

Recall that $e_3$ is a type II edge, $e_1$ and $e_2$ are type I edges of $T$. 
Let $w_i$ be the restriction of $v$ on $e_i\subset \partial T$ as a function of
a single variable, then $w_i$ is a polynomial of degree three which satisfies
\begin{align*}
&w_i^\prime(0)=w_i(0)=w_i(\frac12)=w_i(1)=0\qquad i=1,2,\\
&w^\prime_3(0)=w_3(0)=w_3(1)=w_3^\prime(1)=0.
\end{align*} 
In either case, we conclude $w_i\equiv 0$.

Next, let $z_3$ be the restriction of $\p_{\bar{n}} v$ on $e_3$ as a function 
of a single variable.  Then $z_3$ is a polynomial of degree two satisfying
\begin{align*}
z_3(0)=z_3(\frac12)=z_3(1),
\end{align*} 
which then infers $z_3\equiv 0$.

From the above calculations, we conclude that $(\lambda_3^T)^2$, 
$\lambda_1^T$, and $\lambda_2^T$ are factors of $v$.  However, this is not 
possible, since $v$ is a polynomial of degree three, unless $v\equiv 0$.
The proof is complete.
\end{proof}

Let $V^h_3$ be the finite element space associated with $S_3^h$, that is,
\[
V^h_3=\{ v|_T \in \mathbb{P}_3(T), \, v \mbox{ is continuous at every 
degree of freedom in } \Sigma_T, \, \forall T\in \mathcal{T}_h\}.
\]  
We now show that $V^h_3$ is a subspace of $V$.
 
\begin{thm}\label{ConformingThm1}
There holds the inclusion $V^h_3\subset V$.
\end{thm}

\begin{proof}
Let $v\in V^h_3$. By the proof of Lemma \ref{H2lem}, it suffices to show 
$v$ and $\p_{\bar{n}} v$ are both continuous across interior edges of 
$\mathcal{T}_h$.  Let $T_1$ and $T_2$ be two adjacent triangles with 
common edge $e$, and $w$ be the restriction of $v^{T_1}-v^{T_2}$ along 
$e$ as a function of a single variable.  We then have
\begin{alignat*}{2}
&w^\prime(0)=w_1(0)=w(\frac12)=w(1)=0\qquad &&\text{if $e$ is type I},\\
&w^\prime(0)=w(0)=w(1)=w^\prime(1)=0\qquad &&\text{if $e$ is type II}.
\end{alignat*}
Thus, $w\equiv 0$ and the inclusion $V^h_3\subset C^0(\Omega)\subset H^1(\Omega)$ 
holds.

Next, we observe that if $e$ is a type I edge, then
\begin{align*}
\big[\p_{\bar{n}} v\big]\big|_e = \pm \big[\p_\tau v\big]\big|_{e} =0.
\end{align*} 
Hence, $\p_{\bar{n}} v$ is continuous across $e$. On the other hand, 
if $e$ is a type II edge, let $z$ be the restriction of 
$[\p_{\bar{n}} v]|_e=\p_{\bar{n}} v^{T_1}-\p_{\bar{n}} v^{T_2}$ 
along $e$ as a function of a single variable. Since
\begin{align*}
z(0)=z(\frac12)=z(1),
\end{align*}
and $z$ is a polynomial of degree two, it follow that 
$ \big[\p_{\bar{n}} v\big]\big|_{e}=0$.
So $\p_{\bar{n}} v$ is also continuous across $e$. This then concludes the proof.
\end{proof}

\begin{remark}
We note that $V^h_3\not\subset H^2(\Ome)$ because $V^h_3\not\subset C^1(\Omega)$.
\end{remark}

\subsection{A quartic conforming finite element}

In this subsection, we again assume that $\mathcal{T}_h$ is a 
triangulation of $\Ome$ and every triangle of $\mathcal{T}_h$ 
has {\em two} type I edges. We then define the following quartic 
finite element $S_4^h:=(T,Q_T,\Xi_T)$:
\begin{enumerate}
\item[\rm (i)] $T$ is a triangle with two type I edges,
\item[\rm (ii)] $Q_T=\mathbb{P}_4(T)$, the space of quartic polynomials on $T$,
\item[\rm (iii)]
$\displaystyle{
\Xi_T=\left\{\begin{array}{ll} v(a_i) &1\le i\le 3,\\
v(a_{ii3}),\, v(a_{i33}) & 1\le i \le 2,\\ 
v(b_3),\\
v(a_{123}),\\
\nabla v(a_i)(a_j-a_i)&1\le i\le 2,\ 1\le j\le 3,\ j\neq i,\\
\p_{\bar{n}} v(a_{112}),\, \p_{\bar{n}} v(a_{122}),
\end{array}\right.
}$

where $e_3$ is a type II edge, and $a_{ij\ell}=\frac13\big(a_i+a_j+a_\ell)$.
\end{enumerate}

\begin{figure}[htb]
\begin{center}
\includegraphics[angle=0,width=6.15cm,height=4.75cm]{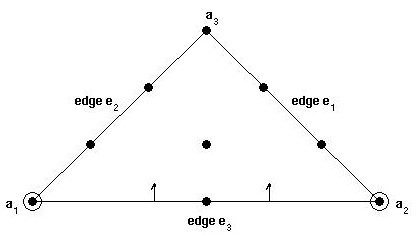}
\end{center}
\caption{\scriptsize Element $S^h_4$. Solid dots indicate function
evaluation, circles indicate first derivative evaluation, and arrows
indicate evaluation of derivatives in the direction $\overline{n}$.}
\end{figure}

\begin{lem}\label{UnisolventLem2}
The set $\Xi_T$ is unisolvent.  That is, any polynomial of degree four 
is uniquely determined by the degrees of freedom in $\Xi_T$.
\end{lem}

\begin{proof}
Suppose $v\in \mathbb{P}_4(T)$ equals zero at all the degrees of 
freedom in $\Xi_T$, and let $w_i$ be the restriction of $v$ to $e_i$ 
as a function of a single variable.  Then
\begin{align*}
&w_i^\prime(0)=w_i(0)=w_i(\frac13)=w_i(\frac23)=w_i(1)=0\qquad i=1,2,\\
&w_3^\prime(0)=w_3(0)=w_3(\frac12)=w_3(1)=w_3^\prime(1)=0.
\end{align*}
Thus, $w_i\equiv 0, i=1,2,3$.

Next, letting $z_3$ be the restriction of $\p_{\bar{n}} v$ on $e_3$
as a function of a single variable, we have
\begin{align*}
z_3(0)=z_3(\frac13)=z_3(\frac23)=z_3(1)=0.
\end{align*}
Hence, $z_3\equiv 0$.

From the above calculations, we conclude that 
$v=a\lambda_1\lambda_2\lambda_3^2$ for some $a\in\mathbf{R}$.  
However, since $0=v(a_{123})=\frac{a}{81}$, we have $a=0$. 
The proof is complete.
\end{proof}

\begin{thm}\label{ConformingThm2}
Let $V^h_4$ be the finite element space associated with $S_4^h$, that is,
\[
V^h_4=\{ v|_T \in \mathbb{P}_4(T), \, v \mbox{ \rm is continuous at every 
degree of freedom in } \Xi_T, \, \forall T\in \mathcal{T}_h\}.
\]
Then there holds the inclusion $V^h_4\subset V$.
\end{thm}

\begin{proof}
Let $v\in V^h_4$, and suppose $T^1,T^2\in\mathcal{T}_h$ are two 
adjacent triangles with common edge $e$.  Let $w$ be the restriction 
of $[v]|_e=v^{T_1}-v^{T_2}$ along $e$ as a function of a single variable, from
\begin{alignat*}{2}
&w^\prime(0)=w(0)=w(\frac13)=w(\frac23)=w(1)=0\qquad &&\text{if $e$ is type I},\\
&w^\prime(0)=w(0)=w(\frac12)=w(1)=w^\prime(1)=0\qquad &&\text{if $e$ is type II},
\end{alignat*}
we conclude $w\equiv 0$. Hence, the inclusion 
$V^h_4\subset C^0(\Omega)\subset H^1(\Omega)$ holds.

If $e$ is a type II edge, we let $z$ denote the restriction of 
$[\p_{\bar{n}} v]|_e=\p_{\bar{n}} v^{T_1}-\p_{\bar{n}} v^{T_2}$ 
along $e$ as a function of one variable. It follows from
\begin{align*}
z(0)=z(\frac13)=z(\frac23)=z(1)=0,
\end{align*}
that $z\equiv 0$.  

Finally, if $e$ is a type I edge, we use the fact that $v$ is continuous 
to conclude
\begin{align*}
\big[\p_{\bar{n}} v\big]\big|_e=\pm \big[\p_\tau v \big]\big|_e=0.
\end{align*}
Thus, $V_4^h\subset V$.
\end{proof}

\begin{remark}
We note that $V^h_4\not\subset H^2(\Ome)$ because $V^h_4\not\subset C^1(\Omega)$. 
\end{remark}

\subsection{Approximation properties of the proposed finite elements}
Let $\Pi^T_k v\in\mathbb{P}_{k}(T)$ denote the standard 
interpolation of $v$ associated with the finite element $S^h_k$,
and define $\Pi^h_k v\in V_k^h$ such that 
$\Pi^h_kv\big|_T=\Pi_k^T(v\big|_T),\ \forall T\in\mathcal{T}_h$.
Before stating the approximation properties of the interpolation 
operator $\Pi_k^T$, we first establish the following technical lemma 
concerning the mesh $\mathcal{T}_h$.
\begin{lem}\label{minanglelem}
Suppose $T\in\mathcal{T}_h$ has two type I edges, and without 
loss of generality, assume $e_3\subset\partial T$ is a type II edge.  
Then there exists a constant $C>0$ that depends only on the minimum 
angle of $T$ such that
\begin{align*}
1-\big(\beta^{(3)}\big)^2\ge C,
\end{align*}
where $\beta^{(3)}=\tau^{(3)}\cdot\overline{n}^{(3)}$.
\end{lem}

\begin{proof}
Since both type I edges of $T$ make an angle of $\frac{\pi}4$
with respect to the $x$-axis (cf. Remark \ref{edgeremark}), 
then there exists $\theta\in (0, \frac{\pi}4]$ such
that the angles of $T$ are $\frac{\pi}2,\ \theta$, and $\frac{\pi}2-\theta$.

Next, we embed $T$ into an isosceles triangle as shown in 
Figure \ref{figmin},
\begin{figure}[htb]
\begin{center}
\includegraphics[scale=0.55]{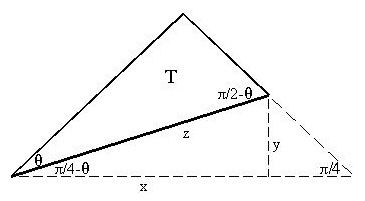}
\end{center}
\caption{\scriptsize Embedding $T$ into an isosceles triangle}
\label{figmin}
\end{figure}
and then obtain 
\begin{align*}
\tau^{(3)}=\frac{(x,y)}{z},\quad \tau^{(3)}\cdot\overline{n}^{(3)}
=\frac{-2xy}{z^2},\quad x=\cos(\frac{\pi}4-\theta)z,\quad 
y=\sin(\frac{\pi}4-\theta)z.
\end{align*}
Hence,
\begin{align*}
\tau^{(3)}\cdot \overline{n}^{(3)}&=\frac{-2xy}{z^2}
=-2\sin(\frac{\pi}4-\theta)\cos(\frac{\pi}4-\theta)
=-\cos(2\theta),
\end{align*}
which implies that
\begin{align*}
1-(\beta^{(3)})^2=1-(\tau^{(3)}\cdot\overline{n}^{(3)})^2
=1-\cos^2(2\theta)=\sin^2(2\theta).
\end{align*}
The proof is complete.
\end{proof}

\begin{remark}
If $\mathcal{T}_h$ is a uniform criss-cross triangulation of $\Omega$, then
$1-(\beta^{(3)})^2=1$ for all type II edges, $e_3$.
\end{remark}

The next theorem establishes the approximation properties of 
the proposed cubic and quartic finite elements.

\begin{thm}\label{affinethm}
For all $m\ge 0$, $p,q\in [1,\infty]$ which are compatible with the inclusion
\begin{align*}
W^{k+1,p}(T)&\hookrightarrow W^{m,q}(T),
\end{align*}
there holds
\begin{align}\label{interpolation}
\|v-\Pi^T_k v\|_{W^{m,q}(T)}
\le C h_T^{k+1-m+\frac{2}{q}-\frac{2}{p}}\|v\|_{W^{k+1,p}(T)}\qquad
\forall v\in W^{k+1,p}(T),
\end{align}
where $h_T={\rm diam}(T)$.
\end{thm}

\begin{proof}
{\em The case $S^h_3$:} Since $S^h_3$ is not an affine family in general, 
the standard scaling technique can not be used directly to prove 
\eqref{interpolation}. To get around this difficulty, the trick is 
to introduce an affine ``relative" of $S^h_3$ and to estimate the 
discrepancy between $S^h_3$ and its ``relative".
To this end, we introduce the following element 
$\mathcal{S}_3^\prime:=(T,P_T,\Sigma_T^\prime)$:
\begin{enumerate}
\item[\rm (i)] $T$ is a triangle with two type I edges, 
\item[\rm (ii)] $P_T=\mathbb{P}_3(T)$,
\item[\rm (iii)]
$\displaystyle{
\Sigma_T^\prime=\left\{\begin{array}{ll} 
             v(a_i) &1\le i\le 3,\\
             v(a_{ik}) & 1\le i \le 2,\\
             \nabla v(a_i)\cdot(a_j-a_i)&1\le i\le 2,\ 1\le j\le 3,\ j\neq i,\\
             \nabla v(b_3)\cdot(a_3-b_3),
\end{array}\right.
}$

where edge $e_3$ is of type II.
\end{enumerate}

\begin{figure}[htb]
\begin{center}
\includegraphics[angle=0,width=5.95cm,height=4.5cm]{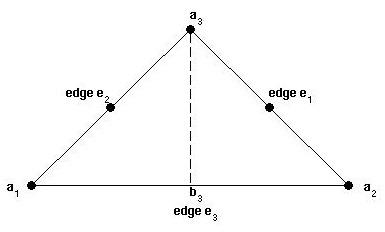}\\
\caption{{\scriptsize Finite element $\mathcal{S}_3^\prime$.}}
\end{center}
\end{figure}

It is easy to see that $\Sigma_T^\prime$ is unisolvent 
in $\mathbb{P}_3(T)$, and that any two triangles are affine equivalent.
Therefore for all $p,q\in [1,\infty]$, $0\le m\le 4$
with $W^{4,p}(T)\hookrightarrow W^{m,q}(T)$, there holds \cite{Ciarlet78}
\begin{align}\label{standard}
\|v-\Lambda^T_3 v\|_{W^{m,q}(T)}
\le C h_T^{4-m+\frac{2}{q}-\frac{2}{p}}\|v\|_{W^{4,p}(T)}
\qquad \forall v\in W^{4,p}(T),
\end{align}
where $\Lambda^T_3$ is the interpolation operator associated with
$\mathcal{S}_3^\prime$.

Define $\Theta^T_3:=\Pi^T_3-\Lambda^T_3$, and note that for 
$v\in W^{4,p}(T)$, $\Theta^T_3v\big|_{e_i}=0$ for $i=1,2,3$. Consequently,
\begin{align*}
&\nabla v(b_3)\cdot (a_3-b_3) = \frac{1}{1-(\beta^{(3)})^2}\Big\{(a_3-b_3)
\cdot \big(\overline{n}^{(3)}-\tau^{(3)}\beta^{(3)}\big)
\p_{\bar{n}^{(3)}}\big(v-\Lambda^T_3v\big)(b_3)\Big\},
\end{align*} 
where $\overline{n}^{(3)}=(n_1^{(3)},-n_2^{(3)})$,
$\beta^{(3)}:=\tau^{(3)}\cdot \overline{n}^{(3)}$,
$n^{(3)}=(n_1^{(3)},n_2^{(3)})$ and $\tau^{(3)}$ denote respectively the unit 
normal and tangential direction of edge $e_3$.

Next, let $q_3$ be the basis function associated with the degree 
of freedom $\nabla v(b_3)(a_3-b_3)$ in $\Sigma_T^\prime$.  We then have
\begin{align*}
\Theta^T_3v=\frac{1}{1-(\beta^{(3)})^2}
\Big\{(a_3-b_3)\cdot \big(\overline{n}^{(3)}-\tau^{(3)}\beta^{(3)}\big)
\p_{\overline{n}^{(3)}}\big(v-\Lambda^T_3v\big)(b_3)\Big\}q_3.
\end{align*}
Therefore, 
\begin{align*}
\|\Theta_3^T v\|_{W^{m,q}(T)} 
\le \frac{1}{1-(\beta^{(3)})^2} \Big\{|a_3-b_3|\cdot
|\overline{n}^{(3)}-\tau^{(3)}\beta^{(3)}|\cdot
\|v-\Lambda_3^Tv\|_{W^{1,\infty}(T)}\|q_3\|_{W^{m,q}(T)}\Big\}.
\end{align*}

Finally, by \eqref{standard} and Lemma \ref{minanglelem} we get
\begin{alignat*}{2}
&1-(\beta^{(3)})^2\ge C,\qquad
&&|a_3-b_3|\le Ch_T,\\
&|\overline{n}^{(3)}-\tau^{(3)}\beta^{(3)}|\le 2,\qquad
&&\|v-\Lambda^T_3v\|_{W^{1,\infty}(T)}\le C h_T^{3-\frac{2}{p}}\|v\|_{W^{4,p}(T)},\\
&\|q_3\|_{W^{m,q}(T)}\le C h_T^{-m+\frac{2}{q}},\qquad &&
\end{alignat*}
where $C$ only depends on the minimum angle of $T$.  Hence,
\begin{align*}
\|\Theta^T_3v\|_{W^{m,q}(T)}
\le C h_T^{4-m+\frac{2}{q}-\frac{2}{p}}\|v\|_{W^{4,p}(T)},
\end{align*} 
and consequently,
\begin{align*}
\|v-\Pi_3^Tv\|_{W^{m,q}(T)} 
\le \|v-\Lambda^T_3v\|_{W^{m,q}(T)}+\|\Theta^T_3v\|_{W^{m,q}(T)}
\le Ch_T^{4-m+\frac{2}{q}-\frac{2}{p}}\|v\|_{W^{4,p}(T)}.
\end{align*}
\smallskip

{\em The case $S_4^h$:} We use a similar argument to show 
\eqref{interpolation} for the element $S_4^h$.  First, we introduce the
following ``relative" $S_4^\prime:=(T,Q_T,\Xi_T^\prime)$ of $S_4^h$:
\begin{enumerate}
\item[\rm (i)] $T$ is a triangle with two type I edges, 
\item[\rm (ii)] $Q_T=\mathbb{P}_4(T)$,
\item[\rm (iii)]
$\displaystyle{
\Xi_T^\prime=\left\{\begin{array}{ll} v(a_i) &1\le i\le 3,\\
           v(a_{ii3}),\, v(a_{ii3}) & 1\le i \le 2,\\
           v(b_3), &\\
           v(a_{123}), &\\
           \nabla v(a_i)\cdot(a_j-a_i)&1\le i\le 2,\ 1\le j\le 3,\ j\neq i,\\
           \nabla v(a_{112})\cdot(a_3-a_{112}),\\
            \nabla v(a_{122})\cdot(a_3-a_{122}), &
\end{array}\right.
}$

where edge $e_3$ is of type II.
\end{enumerate}

\begin{figure}[htb]
\begin{center}
\includegraphics[angle=0,width=5.95cm,height=4.5cm]{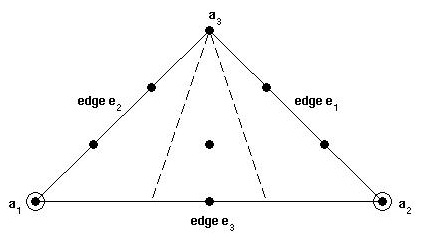}
\caption{\scriptsize Element $\mathcal{S}_4^\prime$. }
\end{center}
\end{figure}

Next, let $\Lambda_4^T$ be the interpolation operator associated with
$\mathcal{S}_4^\prime$, and set $\Theta_4^T:=\Pi^T_4-\Lambda_4^T$.  
Let $r_1$ be the basis function of the element $\mathcal{S}_4^\prime$ 
that is associated with the degree of freedom $\nabla v(a_{112})(a_3-a_{112})$, 
and let $r_2$ be the basis function that is associated with the degree of 
freedom $\nabla v(a_{122})(a_3-a_{122})$. Then for $v\in W^{5,p}(T)$
\begin{align*}
\Theta^T_4v&=\frac{1}{1-(\beta^{(3)})^2}\Big\{(a_3-a_{112})\cdot 
\big(\overline{n}^{(3)}-\tau^{(3)}\beta^{(3)}\big)
\p_{\overline{n}^{(3)}}\big(v-\Lambda^T_4v\big)(a_{112})r_1\\
&\qquad\qquad+ (a_3-a_{122})\cdot \big(\overline{n}^{(3)}-\tau^{(3)}\beta^{(3)}
\big)\p_{\overline{n}^{(3)}}\big(v-\Lambda^T_4v\big)(a_{122})r_2\Big\}.
\end{align*}

Using the fact $\mathcal{S}_4^\prime$ is affine equivalent and 
applying Lemma \ref{minanglelem} we get
\begin{alignat*}{2}
&1-(\beta^{(3)})^2\ge C,\qquad
&&|a_3-a_{112}|,\ |a_3-a_{122}|\le Ch_T,\\
&|\overline{n}^{(3)}-\tau^{(3)}\beta^{(3)}|\le 2,\qquad
&&\|v-\Lambda^T_4v\|_{W^{1,\infty}(T)}
\le C h_T^{4-\frac{2}{p}}\|v\|_{W^{5,p}(T)},\\
&\|r_i\|_{W^{m,q}(T)}\le C h_T^{-m+\frac{2}{q}},\ i=1,2.\qquad &&
\end{alignat*}
Therefore,
\begin{align*}
\|\Theta_4^T\|_{W^{m,q}(T)} 
\le Ch_T^{5-m+\frac{2}{q}-\frac{2}{p}}\|v\|_{W^{5,p}(T)},
\end{align*}
and consequently,
\begin{align*}
\|v-\Pi_4^T\|_{W^{m,q}(T)}
\le \|v-\Lambda_4^Tv\|_{W^{m,q}(T)}+\|\Theta^T_4 v\|_{W^{m,q}(T)}
\le Ch_T^{5-m+\frac{2}{q}-\frac{2}{p}}\|v\|_{W^{5,p}(T)}.
\end{align*}
The proof is complete.
\end{proof}

We note that if a uniform criss-cross mesh is used such that
every triangle has two type I edges (see Figure \ref{squaremesh2}),
then $\nabla v(b_3)(a_3-b_3)=\pm \p_{\bar{n}} v$ in the definition 
of $\Sigma^\prime_T$. This observation leads to the following corollary.
\begin{cor}
Suppose $\mathcal{T}_h$ is the uniform criss-cross triangulation of $\Omega$,
then $S^h_3=S_3^\prime$. Hence, $S^h_3$ is an affine family.  
\end{cor}

\section{Finite element formulation and convergence analysis}\label{sec-3}
Let $V^h_k\, (k=3,4)$ be the finite element subspaces of $V$ constructed 
in the previous section. Define 
\[
V_{k0}^{h}:=\{v\in V^h_k;\ v\big|_{\p\Ome}=\p_{\bar{n}}v\big|_{\p\Ome}=0\}.
\]
Based on the weak formulation \eqref{firstvariational}, we define our 
finite element method for problem \eqref{boxproblem1}--\eqref{boxproblem2}
as seeking $u_h\in V_{k0}^h$ such that
\begin{align} \label{firstconforming}
A^\delta(u_h,v_h)=(f,v_h)\qquad \forall v_h\in V_{k0}^h.
\end{align}

On noting \eqref{coercive}--\eqref{bounded}, an application of 
Cea's Lemma \cite{Ciarlet78} yields the following result.

\begin{lemma}\label{Cea}
There exists a unique solution to \eqref{firstconforming}.
Furthermore, the following error estimate holds:
\begin{align*}
\|u-u_h\|_V\le C \inf_{v_h\in V_{k0}^h} \|u-v_h\|_V.
\end{align*}
\end{lemma}

Combining Lemma \ref{Cea} and Theorem \ref{affinethm} with $p=q=2,\ m=1,2$
we immediately get the following energy norm error estimate.
\begin{thm}\label{Cea2}
If $u\in H^s(\Omega)\ (s\ge 3)$ then
\begin{align*}
\|u-u_h\|_V\le C h^{\ell-2}\big(\sqrt{\delta}+h\big)\|u\|_{H^\ell},
\qquad \ell={\rm min}\{k+1,s\}.
\end{align*} 
\end{thm}

Next, using a duality argument, we obtain an error estimate in the $L^2$-norm.
\begin{thm}\label{L2thm}
Suppose $u\in H^s(\Omega)\ (s\ge 3)$. Then there holds the following 
error estimate:
\begin{align}
\label{C1l2err}\|u-u_h\|_\lt
&\le C\hat{C}_{0,0} h^{\ell-1}\big(\sqrt{\delta}+h\big) \|u\|_{H^\ell}
\qquad \ell={\rm min}\{k+1,s\}.
\end{align}
\end{thm}

\begin{proof}
Denote the error by $e_h:=u-u_h$, and let $\varphi\in V_0$ 
be the solution to the following auxiliary problem:
\begin{align*}
A^\delta(\varphi,v)&=\langle e_h,v\rangle \qquad \forall v\in V_0.
\end{align*}
It follows from Theorems \ref{thm2.1} and  \ref{regularity_thm} 
that the above problem has a unique solution $\varphi$ and 
\begin{align}\label{H2bound}
\sqrt{\delta} \|\nab\Box \varphi\|_{L^2} 
+ \|\Del \varphi\|_{L^2} \leq \hat{C}_{0,0} \|e_h\|_{L^2}.
\end{align}
We then have
\begin{align}\label{e4.4}
\|e_h\|_\lt^2&=A^\delta(e_h,\varphi)
=A^\delta(e_h,\varphi-\mathcal{P}_k^h\varphi)
\leq \|e_h\|_V \|\varphi-\mathcal{P}_k^h\varphi\|_V,
\end{align}
where $\mathcal{P}_k^h$ denotes the $L^2$-projection to $V^h_{k0}$. 

By the definition of $\|\cdot\|_V$ and \eqref{H2bound} we get
\begin{align}\label{e4.5}
\|\varphi-\mathcal{P}_k^h\varphi\|_V 
&\leq \sqrt{\delta} \| \Box \varphi - \mathcal{P}_k^h \Box \varphi\|_{L^2}
 + \| \nab \varphi - \mathcal{P}_k^h \nab \varphi\|_{L^2} \\
&\leq C\sqrt{\delta} h\|\nab \Box \varphi\|_{L^2}  
+ Ch \| \nab(\nab \varphi)\|_{L^2} \nonumber\\
&\leq Ch \bigr(\sqrt{\delta} \|\nab\Box \varphi\|_{L^2}  
+  \|\Del \varphi\|_{L^2} \bigr) \nonumber\\
&\leq C\hat{C}_{0,0} h  \|e_h\|_{L^2}.  \nonumber
\end{align}
Thus, it follows from Theorem \ref{Cea2}, \eqref{e4.4}, and \eqref{e4.5} that
\begin{align*}
\|e_h\|_\lt
\le C\hat{C}_{0,0} h^{\ell-1}(\sqrt{\delta}+h)\|u\|_{H^\ell}.
\end{align*}
The proof is complete.
\end{proof}

We conclude this section with a few remarks.

\begin{remark}
(a) The energy norm error estimate is optimal, on the other hand, the 
$H^1$ and $L^2$ norm estimates are optimal provided that 
$\sqrt{\delta} \simeq h$.


(b) All above convergence results only hold for the restricted 
meshes, that is, every triangle of the mesh $\mathcal{T}_h$ 
needs to have two type I edges.  As already mentioned at the 
end of Section \ref{sec-2.1}, for arbitrary mesh $\mathcal{T}_h$, 
$V^h\subset V$ will implies that $V^h$ (and $V^{h}_0$) needs to be a $C^1$ 
finite element space on $\mathcal{T}_h$ such as Argyris, Hsieh-Clough-Tocher,
Bogner-Fox-Schmit elements (cf. \cite{Ciarlet78}). In such a case, 
it follows from Lemma \ref{Cea} that 
\begin{align*}
\|u-u_h\|_V &\le C \inf_{v_h\in V^{h}_0} \|u-v_h\|_V \\
&\le C\inf_{v_h\in V^{h}_0 }\bigl\{\sqrt{\delta}\|u-v_h\|_{H^2}
+\|u-v_h\|_{H^1} \bigr\} \\
&\le C h^{\ell-2}(\sqrt{\delta}+h)\|u\|_{H^\ell},
\end{align*}
where $\ell={\rm min}\{k+1,s\}$ and $k (\geq 5)$ is the order of the 
$C^1$ finite element.  Thus, we still get optimal order error estimate in the 
energy norm.  Although, as expected, using $C^1$ finite elements is not 
efficient to solve the bi-wave problem (cf. \cite{feng_neilan}).
\end{remark}

\section{Numerical experiments and rates of convergence}

In this section, we provide some numerical experiments to gauge
the efficiency and validate the theoretical error bounds for the 
finite element $S^h_3$ developed in the previous sections. 

\medskip
{\bf Test 1.} For this test, we calculate the rate of convergence 
of $\|u-u_h\|$ for fixed $\delta$ in various norms and compare each 
computed rate with its theoretical estimate. All our computations are done 
on the square domain $\Omega=(0,1)^2$ using the criss-cross mesh. We use the 
source function
\begin{align*}
f(x,y)=&-2048\pi^4\delta\big(\cos^2(4\pi x)-\sin^2(4\pi y)\big)
-32\pi^2\Big\{\sin^2(4\pi y)\big(\cos^2(4\pi x)-\sin^2(4\pi x)\big) \\
&\qquad +\sin^2(4\pi x)\big (\cos^2(4\pi y)-\sin^2(4\pi y)\big)\Big\},
\end{align*} 
so that the exact solution is given by $u(x,y)=\sin^2(4\pi x)\sin^2(4\pi y)$.

\begin{figure}[htb]
\centerline{
\indent\indent\includegraphics[scale=0.14]{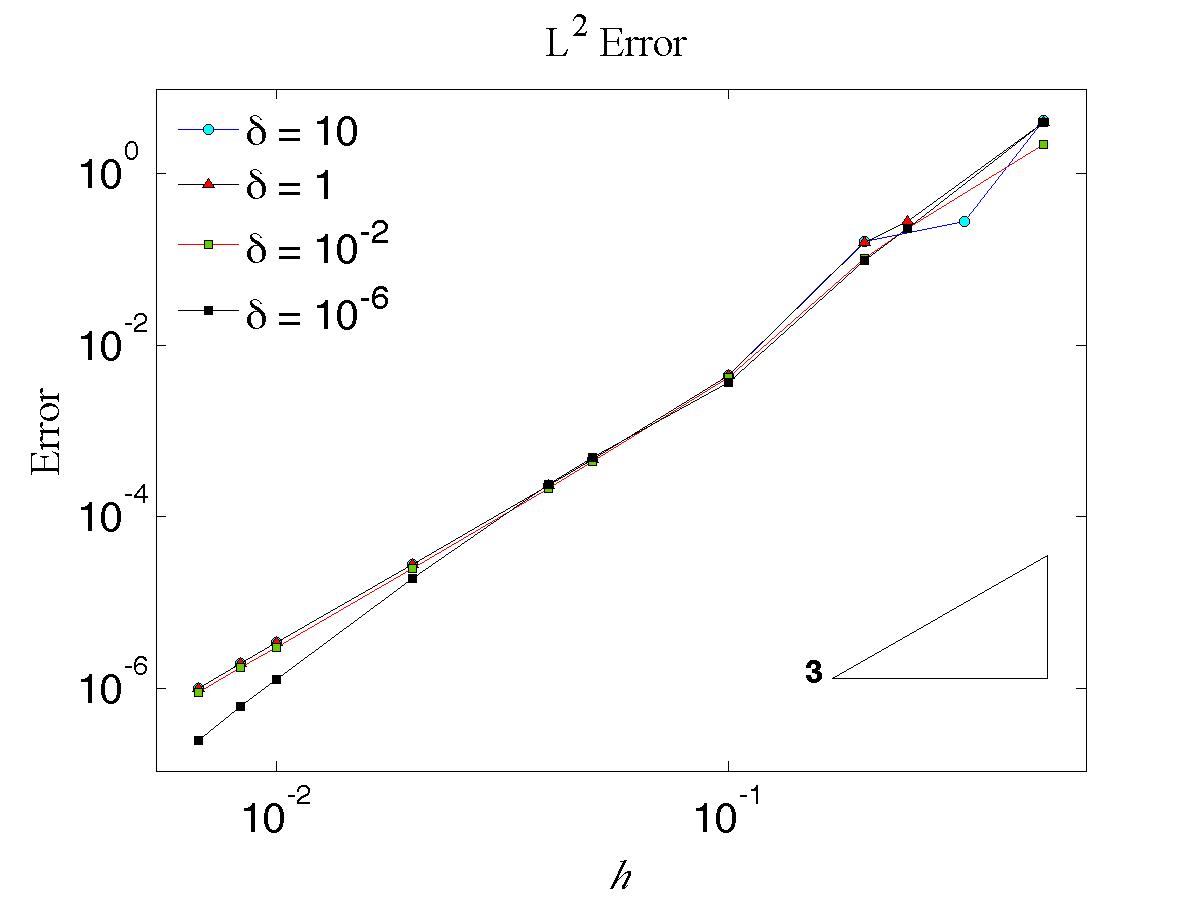}
\indent\indent\includegraphics[scale=0.14]{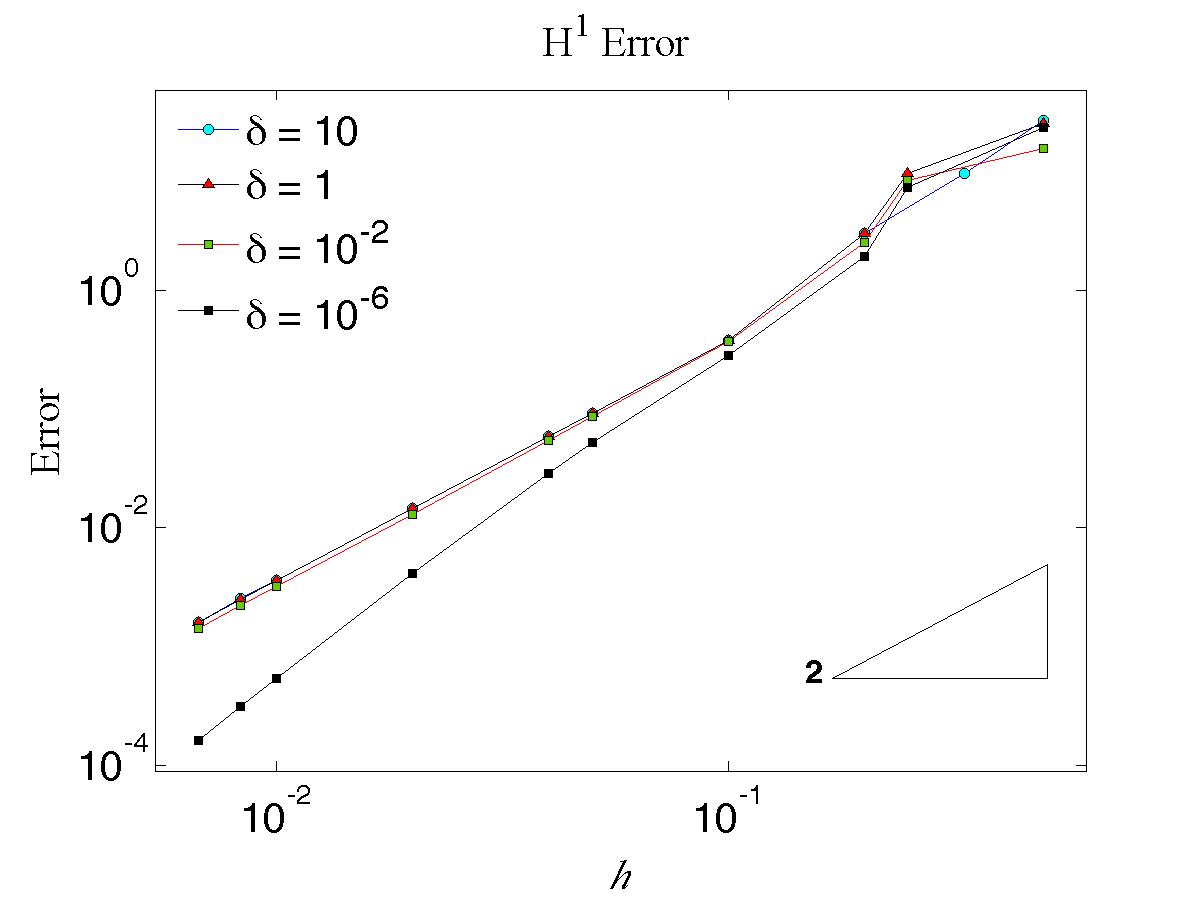}
}
\centerline{
\indent\indent\includegraphics[scale=0.14]{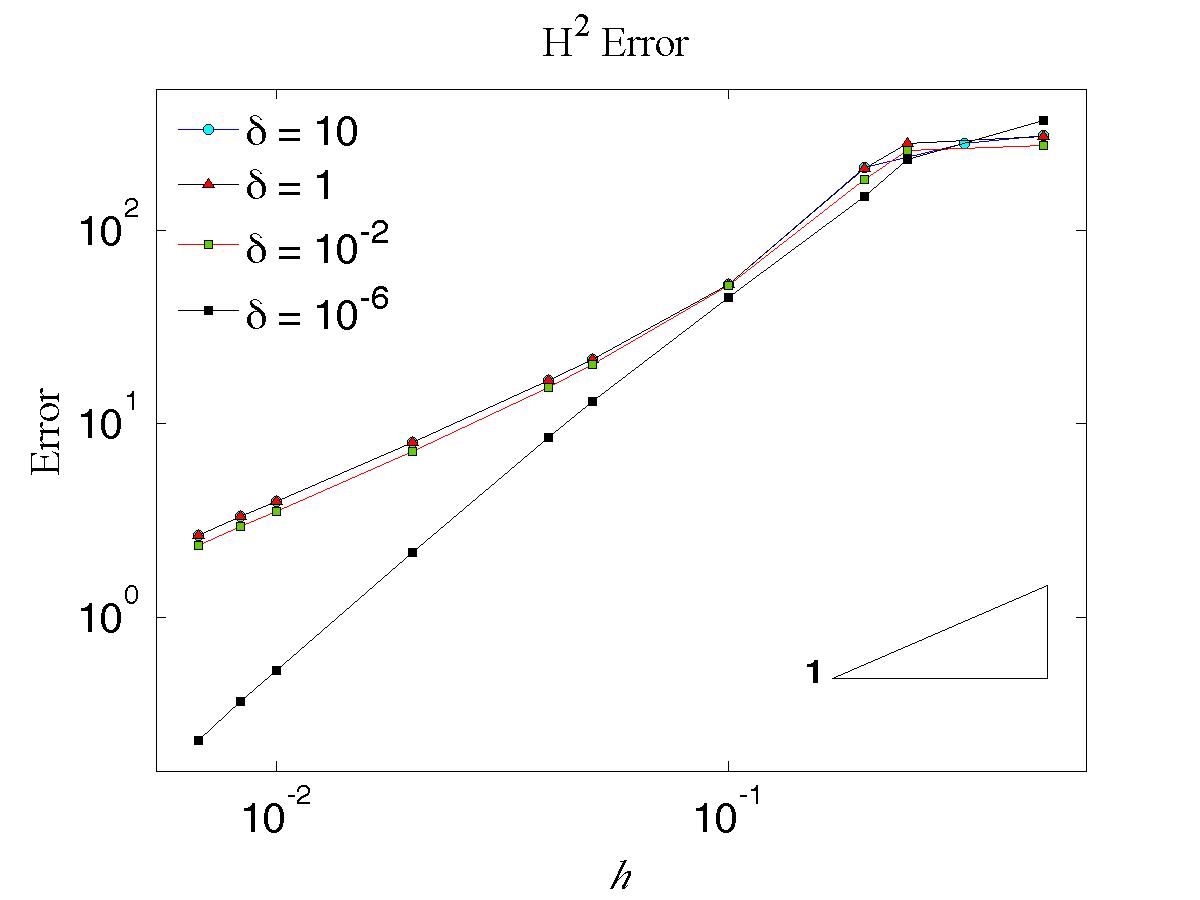}
\indent\indent\includegraphics[scale=0.14]{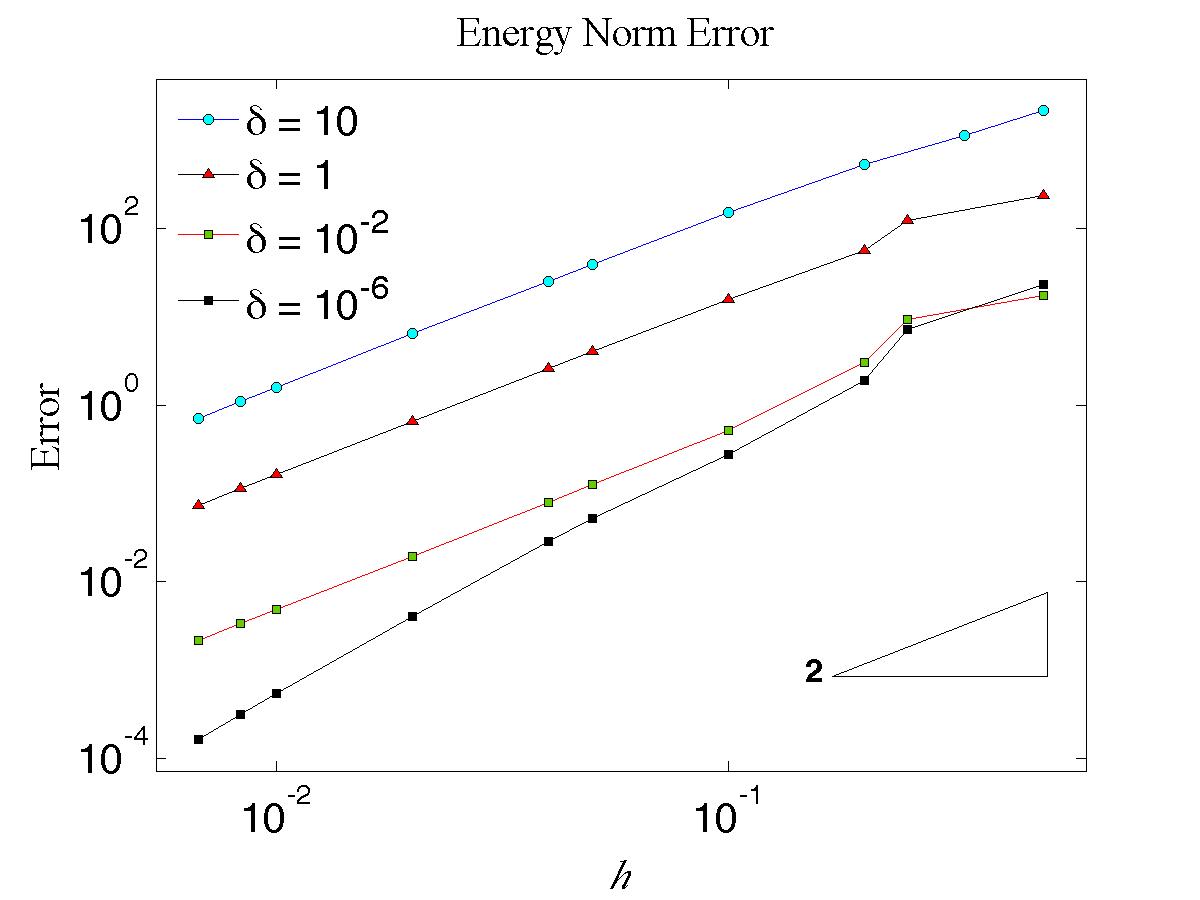}
}
\caption{\label{Test2fig}{\scriptsize Test 1.  $L^2$ norm , $H^1$ norm, 
$H^2$ norm, and energy norm errors with $\delta=10,\ 1,\ 10^{-2}$ and $10^{-6}$.}}
\end{figure}

We list the computed errors in Table \ref{figtest2b} 
for $\delta$-values $10, 1, 10^{-2}$ and $10^{-6}$, and also plot the 
results in Figure \ref{Test2fig}.  As expected, the rates of convergence 
depend on both the parameter $h$ and $\delta$.  In fact, Corollary \ref{Cea2} 
tells us that for $\sqrt{\delta} >> h$
\begin{align*}
\|u-u_h\|_V&\le C h^2(\sqrt{\delta} +h)\|u\|_{H^4}\le Ch^2\|u\|_{H^4},\\
\|u-u_h\|_{H^1}&\le C h^2(\sqrt{\delta} +h)\|u\|_{H^4}\le Ch^2\|u\|_{H^4},\\
\|u-u_h\|_\lt&\le C\hat{C}_{0,0} h^3(\sqrt{\delta} +h)\|u\|_{H^4}
\le C\hat{C}_{0,0} h^3\|u\|_{H^4},\\
\end{align*}
while for $\sqrt{\delta} \leq h$
\begin{align*}
\|u-u_h\|_V&\le C h^2(\sqrt{\delta} +h)\|u\|_{H^4}\le Ch^3\|u\|_{H^4},\\
\|u-u_h\|_{H^1}&\le C h^2(\sqrt{\delta}+h)\|u\|_{H^4}\le Ch^3\|u\|_{H^4},\\
\|u-u_h\|_\lt&\le C\hat{C}_{0,0} h^3(\sqrt{\delta} +h)\|u\|_{H^4}
\le C \hat{C}_{0,0} h^4\|u\|_{H^4}.
\end{align*}
We find that the computed bounds agree with these theoretical bounds.

In addition, although a theoretical proof of the following convergence 
rate has yet to be shown, the computed solutions also indicate that
\begin{align*}
\|u-u_h\|_{2,h}&\le Ch(\sqrt{\delta} +h)\|u\|_{H^4},
\end{align*}
where 
\[
\|u-u_h\|_{2,h}^2:=\sum_{T\in\mathcal{T}_h} \|u-u_h\|_{H^2(T)}^2.
\]

\begin{figure}[htb]
\centerline{
\includegraphics[scale=0.25]{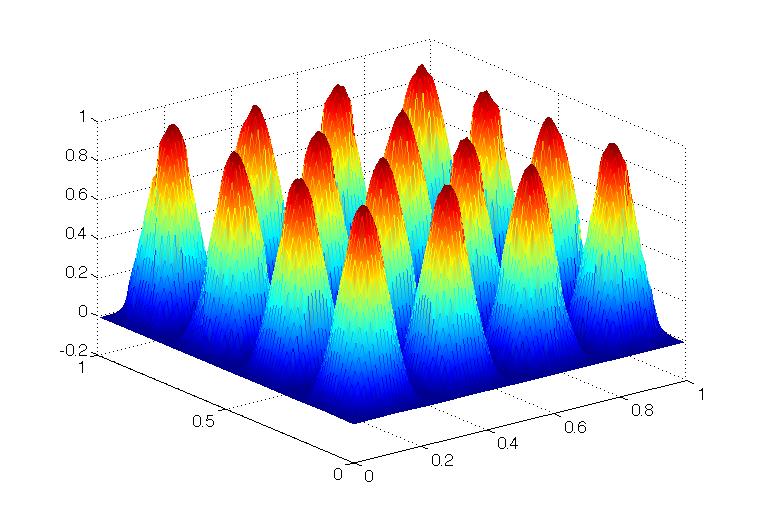}
\includegraphics[scale=0.25]{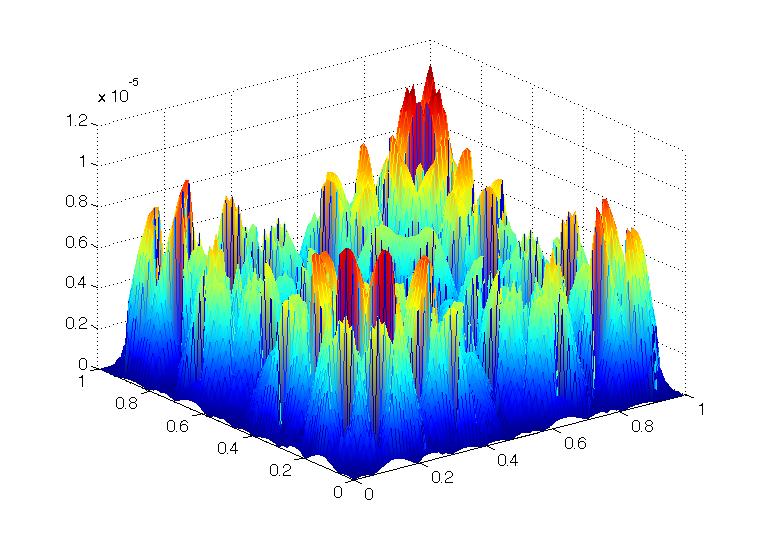}
}
\caption{\label{Compfig}{\scriptsize Test 1.  Computed solution (left) and error (right)
 with $\delta=10^{-2}$ and $h=0.01$.}}
\end{figure}

{\small
\begin{table}[htb]
\begin{center}{\scriptsize
\begin{tabular}{rrrrrr}
$\delta$ &\multicolumn{1}{c}{$h$} & $\|\cdot\|_{L^2}$ err. (cnv. rate) 
&$\|\cdot\|_{H^1}$ err.(cnv. rate) & $\|\cdot\|_{h,2}$ err. (cnv. rate)
&$\|\cdot\|_{V}$ err.(cnv. rate)\\
\noalign{\smallskip}\hline\noalign{\smallskip}
10 &0.5000&4.17($-$)&26.4($-$)& 311.62($-$) &2191.62($-$)\\
&0.3333	& 2.76E-01(6.694)&	9.54(2.514)&	284.05(0.228)& 1147.49(1.596)\\
&0.2000	&1.59E-01(1.079)&	2.99(2.273)&	 211.50(0.577) &535.69(1.491)\\
&0.1000&4.41E-03(5.176)	&3.75E-01(2.995)&	52.54(2.009) &153.70(1.801)\\
&0.0500	&4.64E-04(3.248)&	9.11E-02(2.041)&	21.61(1.282) & 39.84(1.948)\\
&0.0400	&2.31E-04(3.117)&5.82E-02(2.010)&	16.79 (1.130) & 25.61(1.980)\\
&0.0200	&2.79E-05(3.054)&1.45E-02(2.004)& 8.05(1.060) &6.44(1.992)\\
&0.0100&3.45E-06(3.014)	&3.62E-03(2.001)	&3.98(1.016)&1.61(1.998)\\
&0.0083	&1.99E-06(3.006)&	2.52E-03(2.000)& 3.31(1.006)&1.12(1.999)\\
&0.0067&	1.02E-06(3.004)&1.61E-03(2.000)	&2.65(1.004)&0.72(1.999)\\
\noalign{\smallskip}\hline\noalign{\smallskip} 1 &
0.5000	&3.93($-$)	&25.3($-$)	&306.88($-$)&238.43($-$)\\
&0.2500	&2.75E-01(3.837)&9.52(1.413)&	283.58(0.114)& 123.22(0.952)\\
&0.2000	&1.57E-01(2.523)&2.98(5.210)&210.95(1.326)	&56.24(3.515)\\
&0.1000	&4.40E-03(5.152)&3.75E-01(2.989)&52.53(2.006)	&15.71(1.840)\\
&0.0500	&4.63E-04(3.249)&9.09E-02(2.043)&21.58(1.284)	&4.07(1.950)\\
&0.0400	&2.31E-04(3.118)&5.81E-02(2.012)&16.76(1.132)	&2.61(1.981)\\
&0.0200	&2.78E-05(3.055)&1.45E-02(2.005)&8.03(1.061)	&0.66(1.992)\\
&0.0100	&3.44E-06(3.015)&3.61E-03(2.001)&3.97(1.016)	&0.16(1.998)\\
&0.0083	&1.99E-06(3.006)&2.51E-03(2.000)&3.31(1.006)	&0.11(1.999)\\
&0.0067	&1.02E-06(3.004)&1.61E-03(2.000)&2.64(1.004)	&0.07(1.999)\\
&0.0056	&5.89E-07(2.999)&1.12E-03(2.000)&2.20(1.003)	&0.05(1.998)\\
\noalign{\smallskip}\hline\noalign{\smallskip}
$10^{-2}$ &
0.5000 & 2.15($-$)	&15.4($-$)	&276.56($-$)&17.4($-$)\\	
&0.3333 & 2.25E-01(3.259)&	8.38(0.879)&	260.32(0.087)&	9.48(0.877)\\
&0.2000 & 1.02E-01(3.556)&	2.53(5.365)&	183.36(1.571)&	3.07(5.047)\\
&0.1000 &4.21E-03(4.597)&3.69E-01(2.780)&52.08(1.816)&5.22E-01(2.558)\\
&0.0500 &4.36E-04(3.269)&8.55E-02(2.107)&20.31(1.358)&1.25E-01(2.058)\\
&0.0400 & 2.15E-04(3.175)&5.38E-02(2.082)&15.52(1.207)&7.93E-02(2.049)\\
&0.0200 & 2.50E-05(3.101)&1.30E-02(2.049)&7.21(1.106)&1.94E-02(2.030)\\
&0.0100 & 3.06E-06(3.033)&3.21E-03(2.016)&3.53(1.030)&4.82E-03(2.010)\\
&0.0083 & 1.77E-06(3.013)&2.23E-03(2.006)&2.94(1.012)&3.35E-03(2.004)\\
&0.0067 & 9.02E-07(3.009)&1.42E-03(2.005)&2.34(1.008)&2.14E-03(2.003)\\
\noalign{\smallskip}\hline\noalign{\smallskip}
$10^{-6}$ &
0.5000	&3.93($-$)		&23.3($-$)	&	374.18($-$)			&23.3($-$)	\\
&0.2500	&2.28E-01(4.108)&	7.25(1.686)&233.21(0.682)&7.25(1.686)\\
&0.2000	&9.68E-02(3.831)&	1.93(5.929)&	149.75(1.985)	&1.93(5.929)\\
&0.1000	&3.70E-03(4.708)&	2.81E-01(2.782)&45.13(1.731)&2.81E-01(2.782)\\
&0.0500	&4.92E-04(2.914)&	5.21E-02(2.429)&13.09(1.786)&5.21E-02(2.429)\\
&0.0400	&2.35E-04(3.298)&	2.89E-02(2.647)&8.54(1.915)&2.89E-02(2.647)\\
&0.0200	&1.91E-05(3.626)&	4.11E-03(2.813)&2.15(1.989)&4.11E-03(2.813)\\
&0.0100	&1.28E-06(3.898)&	5.39E-04(2.931)&0.53(2.024)&5.39E-04(2.931)\\
&0.0083	&6.19E-07(3.981)&	3.15E-04(2.943)&0.36(2.035)&3.15E-04(2.943)\\
&0.0067	&2.53E-07(4.005)&	1.64E-04(2.934)&0.23(2.039)&1.64E-04(2.933)
\end{tabular}
}
\caption{\label{figtest2b}\scriptsize Test 1.  Errors 
with estimated rates of convergence}
\end{center}
\end{table}
}

\begin{figure}[htb]
\centerline{
\indent\indent\includegraphics[scale=0.25]{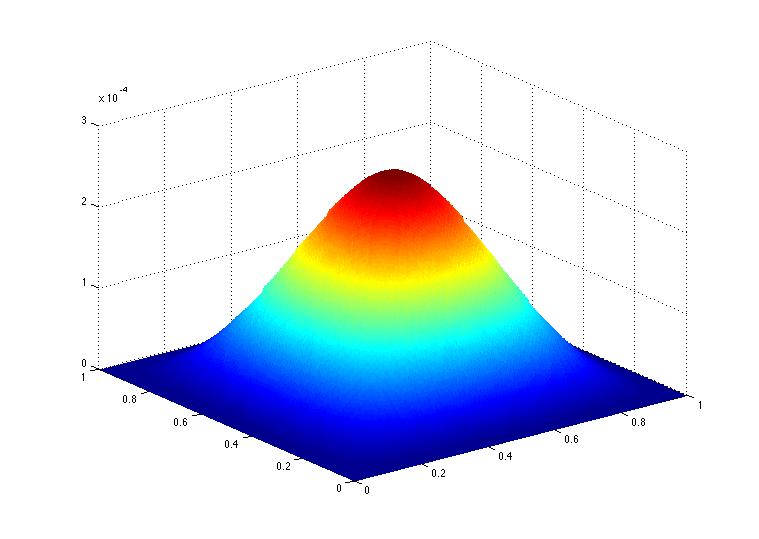}
\indent\indent\includegraphics[scale=0.25]{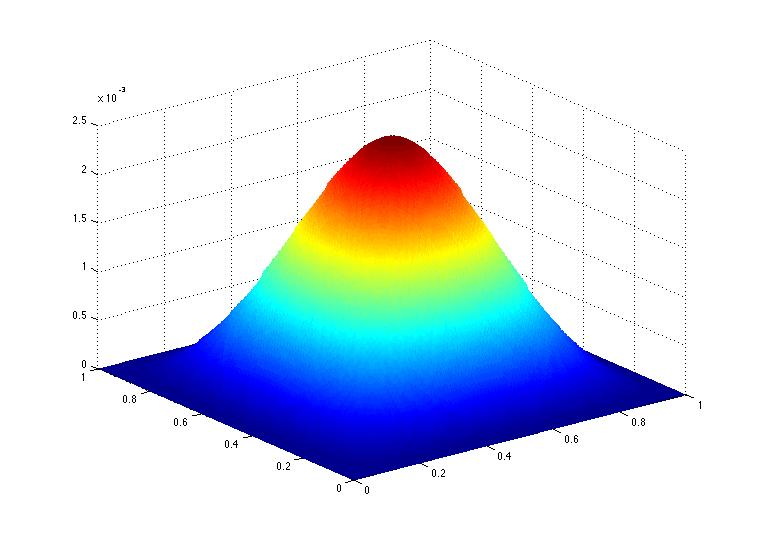}
}
\centerline{
\indent\indent\includegraphics[scale=0.25]{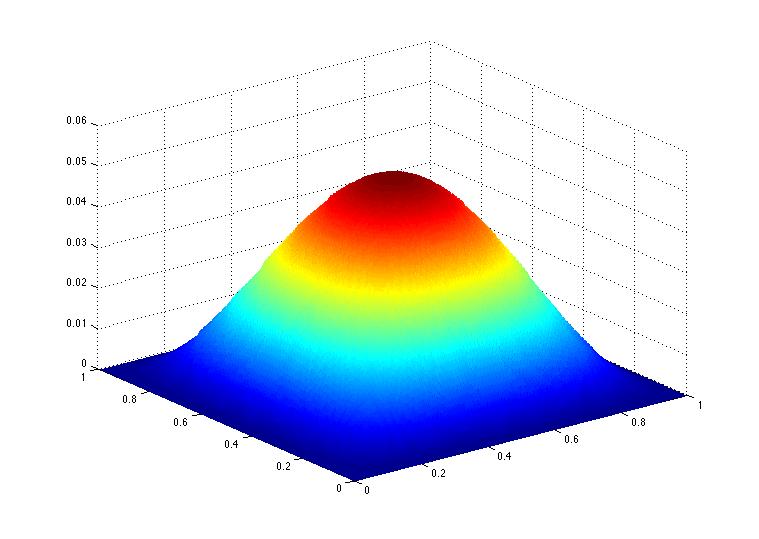}
\indent\indent\includegraphics[scale=0.25]{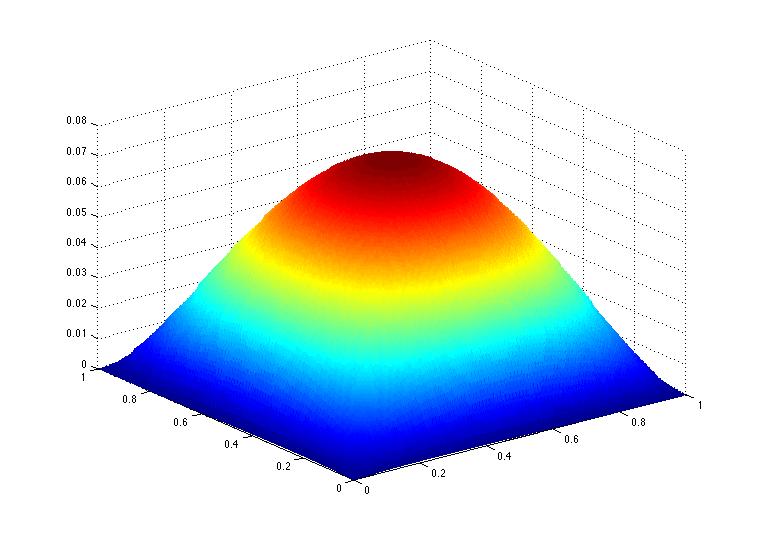}
}
\caption{\label{Test2d}{\scriptsize Test 2.  Computed solution 
with source function $f=1$ and $h=0.01$ with $\delta=10$(top left), 
$\delta=1$(top right), $\delta=10^{-2}$(bottom left), and 
$\delta=10^{-6}$(bottom right).}}
\end{figure}

\medskip
{\bf Test 2.}
This test is the same as the first, but we now use the following source function:
\begin{align*}
f=1.
\end{align*} 
We note that the exact solution is unknown. We plot the solution 
with $h=0.01$ and $\delta$-values $10, 1, 10^{-2},$ and $10^{-6}$ 
in Figure \ref{Test2d}. As expected, the solution is more and more like
the solution of the corresponding Poisson problem as $\delta$ gets smaller
and smaller.

\bigskip
\noindent {\bf Acknowledgments.}
The work of both authors was partially supported by the NSF grant DMS-0710831.
The authors would like to thank Professor Qiang Du of Penn State University
for bringing the bi-wave problem to their attention and for providing 
the relevant references on $d$-wave superconductors.
 

\end{document}